\documentclass[a4paper,11pt,mathcal]{article}
\usepackage{amsfonts,eucal,amssymb,amsmath,amsthm,amscd} 
\newtheorem{theorem}{Theorem}[section]
\newtheorem{corollaire}[theorem]{Corollary}
\newtheorem{remarque}[theorem]{Remark}
\newtheorem{proposition}[theorem]{Proposition}
\newtheorem{lemme}[theorem]{Lemma}
\newtheorem{definition}[theorem]{Definition}

\newenvironment{preuve}{\medskip \noindent {\bf Proof: }}
   {$\diamondsuit$ }

\newcommand{\Conf}{\mbox{Conf }}
\newcommand{\Is}{\mbox{Is }}

\newcommand{\Ad}{\text{Ad }}
\newcommand{\Aut}{\text{Aut }}

\newcommand{\lieo}{\ensuremath{\mathfrak{o}}}

\newcommand{\liep}{\ensuremath{\mathfrak{p}}}

\newcommand{\lieg}{\ensuremath{\mathfrak{g}}}
\newcommand{\lien}{\ensuremath{\mathfrak{n}}}

\newcommand{\re}{\ensuremath{{\bf R}}}
\newcommand{\CC}{\ensuremath{{\bf C}}}
\newcommand{\NN}{\ensuremath{{\bf N}}}
\newcommand{\XX}{\ensuremath{{\bf X}}}

\newcommand{\X}{\ensuremath{{\bf X}}}

\newcommand{\ho}{\ensuremath{{ \hat \Omega}}}

\newcommand{\hm}{\ensuremath{{ \hat M}}}
\newcommand{\hL}{\ensuremath{{ \hat L}}}
\newcommand{\hn}{\ensuremath{{ \hat N}}}
\newcommand{\hx}{\ensuremath{{ \hat x}}}
\newcommand{\hy}{\ensuremath{{ \hat y}}}
\newcommand{\hz}{\ensuremath{{ \hat z}}}

\newcommand{\ol}{\ensuremath{{ \omega^L}}}
\newcommand{\om}{\ensuremath{{ \omega^M}}}
\newcommand{\on}{\ensuremath{{ \omega^N}}}

\newcommand{\oo}{\ensuremath{{ \omega^{\Omega}}}}

\newcommand{\pcho}{\ensuremath{{ \partial_c{ \hat \Omega}}}}
\newcommand{\pho}{\ensuremath{{ \partial{ \hat \Omega}}}}
\newcommand{\psm}{\ensuremath{{ \partial_s{ M}}}}
\newcommand{\pshm}{\ensuremath{{ \partial_{\sigma} {\hat M}}}}
\newcommand{\phm}{\ensuremath{{ \partial {\hat M}}}}
\newcommand{\pchm}{\ensuremath{{ \partial_c{\hat M}}}}

\newcommand{\hl}{\ensuremath{{ {\hat \Lambda}}}}
\newcommand{\hlc}{\ensuremath{{ {\hat \Lambda}_c}}}

\newcommand{\gga}{\ensuremath{{ G_{\Gamma}}}}
\newcommand{\hlg}{\ensuremath{{ {\hat L}_{\Gamma}}}}
\newcommand{\Lg}{\ensuremath{{ L_{\Gamma}}}}

\begin{document}
\pagenumbering{arabic}
\title{Rigidity at the boundary for conformal structures and other  Cartan geometries}
\author{Charles Frances}
\date{}
\maketitle

\section{Introduction}

Building a conformal boundary is a useful tool to study the conformal properties at infinity of noncompact pseudo-Riemannian manifolds. Mostly motivated by General Relativity, some constructions were proposed by Geroch, Kronheimer and Penrose in \cite{gkp}, and also by Schmidt in \cite{schmidt2}, which attach an abstract conformal boundary to a Lorentzian manifold (recall that a metric on a $n$-dimensional manifold is called Lorentzian when its signature is $(1,n-1)$, i.e $(-+...+)$). The construction in \cite{gkp} uses the causal properties of the spacetime under consideration, so that it does not generalize to other structures than Lorentz ones, while the {\it $b$-boundary} of \cite{schmidt2} can be defined for a lot of structures.  Both constructions are intrinsic, and present the advantage of being quite general.  Their main drawback is that the  boundaries obtained in this way   often present topological pathologies (they may not be Hausdorff).  Also, it is generally very hard to determine those boundaries explicitely.  

Another approach to the conformal boundary problem is to consider a type-$(p,q)$ pseudo-Riemannian manifold $(M,g)$, and to embed it strictly  into another one $(N,h)$, of same signature (and in particular of same dimension) thanks to an embedding $s : (M,g) \to (N,h)$ which is {\it conformal}.  This means that $s^*h=e^{\sigma}g$ for some smooth $\sigma : M \to \re$ (we are dealing with smooth structures in all the paper).  If such an embedding $s$ exists, we can consider $\partial_sM$, the topological boundary of $s(M)$ in $N$ (which is nonempty since $s$ is strict), as a conformal boundary for $(M,g)$.  
A lot of classical objects in Riemannian or pseudo-Riemannian geometry appear naturally endowed with a conformal boundary thanks to such embeddings: for example the $n$-dimensional real hyperbolic space is conformally embedded as the upper-hemisphere of the round sphere ${\bf S}^n$, what attaches to him naturally a $(n-1)$-dimensional sphere as conformal boundary.

This point of view raises two natural questions, which will be the main themes of this article. The first question is: {\it given $(M,g)$ a type-$(p,q)$ pseudo-Riemannian manifold, when does it exist a strict conformal embedding $s: (M,g) \to (N,h)$ (strict means $s(M) \not = N$) into a type-$(p,q)$ manifold $(N,h)$?}  This question leads to the natural notion of {\it conformally maximal manifolds}, namely pseudo-Riemannanian manifolds $(M,g)$ such that any conformal embedding $s : (M,g) \to (N,h)$ (with $\text{dim} M= \text{dim} N$) is onto.  Among the main results of this paper , we will exhibit quite a wide class of Riemannian or pseudo-Riemannian manifolds which are conformally maximal.

The second question which will interest us is to determine to what extent the extrinsic construction of a boundary by conformal embedding (when such embeddings exist) is actually intrinsic.  To make  the problem more precise, let us consider a pseudo-Riemannian manifold $(M,g)$ which is not conformally maximal, and assume that we have two strict conformal embeddings  $s_1: (M,g) \to (N_1,h_1)$ and $s_2 : (M,g) \to (N_2,h_2)$ (let us recall that we always assume $\text{dim} M= \text{dim} N_1= \text{dim} N_2$). We will adress the following question: can the topological boundaries  $\partial_{s_1}M$  and $\partial_{s_2}M$  be very different?
 For example, if $\partial_{s_1}M$ is a nice, smooth submanifold of $N_1$, could it happen that $\partial_{s_2}M$
 is very wild?   If one looks at the 2-dimensional case, the answer seems to be definitely yes. Indeed, the Riemann mapping theorem ensures the existence of conformal embeddings of the Poincar\'e disc into $\CC$, with  very different boundaries.   In constrast with this situation,  our aim in this paper is to show that  rigidity phenomena appear for conformal embeddings of Riemannian manifolds  as soon as the  dimension is at least $3$.   Let us mention that many works already exist about this kind of problem.  For the existence of conformal embedding, we should quote \cite{herzlich}, \cite{carron-herzlich}, \cite{carron-herzlich3}, \cite{schoen-yau} among others, which give sufficient conditions (for example on the Ricci curvature of $g$)  ensuring the existence of  a strict conformal embedding 
$s : (M,g) \to (N,h)$.  Also, in  \cite{anderson1}, \cite{anderson2}, \cite{bibi}, \cite{chrusciel}, \cite{graham-lee}, \cite{klo}, \cite{klo2}, related questions of existence and unicity of the conformal boundary are studied.

Both of the two questions above make sense not only for conformal structures, but for general geometric structures.  As the title suggests it, we will consider them in the quite wide framework of Cartan geometries.  We will define precisely Cartan geometries a little bit later.  For the moment, let us just mention that  they include pseudo-Riemannian metrics, conformal pseudo-Riemannian structures, $CR$-structures, affine, projective structures etc....  We now detail the main results of the paper, begining with the conformal Riemannian framework, which is maybe the most familiar to the reader, and then giving the general results about emebddings of Cartan geometries.


\subsection{Conformally maximal manifolds}

In the metric context, it is quite clear that a complete pseudo-Riemannian manifold is maximal in the sense that any isometric embedding has to be onto.  Of course, the example of Euclidean space shows that completeness is no longer sufficient to imply conformal maximality, which makes this latter notion harder to understand. Thus, exhibiting large classes of conformally maximal structures seems  already interesting.  In this paper, we will prove:  

\begin{theorem}
\label{courbure-constante-maximale}
\begin{enumerate}
\item{Every complete flat Riemannian manifold of dimension $n \geq 3$ which is not conformally equivalent to the Euclidean space is conformally maximal. }
\item{For a complete hyperbolic manifold $M=\Gamma \backslash {\bf H}^n$, $n \geq 3$, the following conditions are equivalent:
\begin{enumerate}
\item{$M$ is conformally maximal.}
\item{$M$ is projectively maximal.}
\item{The limit set $\Lambda_{\Gamma}$ is equal to ${\bf S}^{n-1}=\partial {\bf H}^n$.}
\end{enumerate}}
\end{enumerate}

\end{theorem}

We refer to \cite{kapovich} for the definition of the limit set of a subgroup of $\Is {\bf H}^n$.
 Let us  remark that theorem \ref{courbure-constante-maximale} implies that complete hyperbolic manifolds of finite volume and dimension $\geq 3$ are conformally maximal.  This is in sharp contrast with the $2$-dimensional case where it is possible to ``fill" the cusps, so that conformally, a hyperbolic surface of finite volume is always obtained by removing points from a compact Riemann surface. 

Theorem \ref{courbure-constante-maximale} is actually a particular case  of Theorem \ref{klein-maximal}, which is much more general in the sense that it does not only deal with Riemannian conformal structures, but  with Cartan geometries in general (see section \ref{introduction-maximal} and \ref{applications}). 

Before stating the next result, which is more specific to the Riemannian case, let us define a {\it conformally homogeneous manifold} as a Riemannian manifold $(M,g)$, the conformal group of which acts transitively.

\begin{theorem}
\label{homogene}
Let $(M,g)$ be a  Riemannian manifold of dimension $n \geq 3$, which is conformally homogeneous.  Then it is  conformally   maximal except in the  following cases:

\ \ $1.$ $(M,g)$ is conformally equivalent to the Euclidean space $\re^n$.

\ \ $2.$  $(M,g)$ is conformally equivalent to the real hyperbolic space ${\bf H}^n$.

\ \ $3. $ $(M,g)$ is conformally equivalent to the product ${\bf H}^m \times {\bf S}^k$, with $m \geq 1$, $k \geq 1$, $m+k \geq 3$ (with the convention that ${\bf H}^1$ is the $1$-dimensional Euclidean space).

\ \ $4. $ $(M,g)$ is conformally equivalent to the product ${\bf H}^m \times \re$, $m \geq 2$.
\end{theorem}

Actually, we don't know if the last example can be removed from the list, i.e if ${\bf H}^m \times \re$ is conformally maximal (we gess it is the case).  The three first examples are conformally equivalent to open subsets of the sphere, hence are not conformally maximal.  Let us mention two consequences of the previous result.  A left invariant Riemannian metric on a connected Lie group $G$ is always conformally maximal, exept when $G=\re^n$.  Also, all Riemannian symmetric spaces but the Euclidean and the hyperbolic spaces are conformally maximal.

When a type-$(p,q)$ pseudo-Riemannian manifold admits a strict conformal embedding $s: (M,g) \to (N,h)$, with $(N,h)$ also of type-$(p,q)$, and if $s(M)$ has compact closure in $(N,h)$, we say that $s$ yields a {\it conformal compactification } of $(M,g)$. Let us now quote a  result about existence of conformal compactification:

\begin{theorem}
\label{gros-isometries}
Let $(M,g)$ be a connected Riemannian manifold of dimension $n \geq 3$. 
Assume that the isometry group of $(M,g)$ is noncompact, and that $(M,g)$ is not conformally flat. Then, there is no conformal compactification for $(M,g)$. 
\end{theorem}

Thus, for generic Riemannian manifolds $M$, a product $\re \times M$ does not admit any conformal compactification.

\subsection{Rigidity of conformal boundaries}
\label{introduction-rigidity}

We now discuss the second topic of this article, which is unicity, or at least ``rigidity"  of the conformal boundary defined by a conformal embedding.  As we already mentioned it, the situation in dimension $\geq 3$ becomes much more rigid than it is dor surfaces. As an example of this rigidity, we will prove:


\begin{theorem}
\label{codimension2}
Let $(L,g)$ be a compact Riemannian manifold of dimension $n \geq 3$. Let $M \subsetneq L$ be a strict open subset. We assume that the Hausdorff dimension of $\partial M$ is strictly less than $n-1$. Then if $s : (M,g) \to (N,h)$ is a conformal embedding, there is $M^{\prime}$ an open subset of $L$ containing $M$, and $s^{\prime} : (M^{\prime},g) \to (N,h)$ a conformal diffeomorphism which extends $s$. In particular, if $N$ is compact, then $(N,h)$ is conformally diffeomorphic to $(L,g)$.
\end{theorem}

So, basically, this theorem classify  all possible conformal embeddings of the manifold $(M,g)$ into another Riemannian manifold of the same dimension.  Such a strong conclusion is possible because of the condition on the Hausdorff dimension of the boundary $\partial M$.  In particular, theorem \ref{codimension2} does not deal with domains having a hypersurface as boundary, a case on which we focus now.  

In the forthcoming statements,  we will consider a smooth Riemannian manifold $(L,g)$, and $M \subsetneq L$ an open subset satisfying the condition:\\
\ \\
{ $\bf H_1$}. {\it The closure $\overline M$  of $M $ in $L$ is a compact topological submanifold with  boundary, and the boundary $\partial M$ is locally the graph of a Lipschitz function. }\\
\ \\
By topological submanifold with boundary, we mean here that there is a ${\cal C}^0$ atlas $(U_i,\phi_i)_{i \in I}$ of $L$ such that $\phi_i (U_i) = \re^n$ and if $U_i \cap M \not = \emptyset$, $\phi_i(U_i \cap M)$ is the open hupper half-space ${x_n>0}$.  

The condition on the boundary is that there is a {\it smooth} atlas $(V_i,\psi_i)$ of $L$ such that in each chart $V_i \cong \re^{n-1} \times \re$ intersecting $\partial M$, $\partial M \cap V_i$ is the graph of a Lipschitz function $f_i : \re^{n-1} \to \re$. 

We will say that the boundary $\partial M$ is of class ${\cal C}^{k,\alpha}$, $k \in {\NN}$, $\alpha \in [0,1]$, if the functions $f_i$ are $k$-times differentiable with a $\alpha$-H\"older $k$-th differential.  Note that a ${\cal C}^{0,1}$-function is just a Lipschitz function.  In the article, by a ${\cal C}^{k,\alpha}$-hypersurface of a manifold, we will denote a subset which is locally the graph of a ${\cal C}^{k,\alpha}$-function. 

We can now state our first theorem:

\begin{theorem}
\label{immersion}
Let $(L,g)$ be a smooth Riemannian manifold of dimension $n \geq 3$, and $M \subsetneq L$ an open subset  satisfying the hypothesis $H_1$ above. We  assume that the boundary $\partial M$ is  of class ${\cal C}^{k,\alpha}$, $k + \alpha\geq 1$. Let  $s : (M,g) \to (N,h)$ be a strict conformal embedding. 
\begin{enumerate}
\item{There is a nonempty open set $\Lambda \subset \partial M$ such that $s$ extends to a ${\cal C}^{k,\alpha}$- immersion:
$$\partial s: \Lambda \to N  $$ whose image is $\partial_sM$.  When $k \geq 1$,  $\partial s$ is a conformal map and its fibers have at most two elements.}
\item{If $s(M)$ has compact closure in $N$, then $\Lambda = \partial M$.  There is a dense open set $\partial_s^{o}M \subset \partial_sM$ which is a ${\cal C}^{k,\alpha}$-hypersurface of $N$, and a dense open subset $\Lambda^{o} \subset \partial M$, such that $\partial s (\Lambda^o)=\partial_s^{o}M $.}
\end{enumerate}
\end{theorem}

This result can be seen as a generalization to dimension $\geq 3$ of Carath\'eodory's theorem about the boundary behaviour of conformal mappings (see \cite{caratheodory}). In Theorem \ref{immersion}, no assumption is made on the image of the embedding. If we add some (quite mild) assumptions on this image, the result can be precised:

\begin{theorem}
\label{geometrique}
Let $(L,g)$ be a smooth Riemannian manifold of dimension $n \geq 3$, and $M \subsetneq L$ an open subset  satisfying the hypothesis $H_1$ above. We assume that the boundary $\partial M$ is of class ${\cal C}^{k,\alpha}$, $k+\alpha \geq 1$. 
Let  $s : (M,g) \to (N,h)$ be a strict conformal embedding. If $\partial_sM$ is locally the graph of a ${\cal C}^0$ function, then:
\begin{enumerate}
\item{ $\partial_sM$ is a ${\cal C}^{k,\alpha}$-hypersurface of $N$.}
\item{If $s(M)$ has compact closure in $N$, then $s$ extends to  a ${\cal C}^{k,\alpha}$-immersion  $\partial s : \partial M \to \partial_s M$ which is onto. In restriction to each connected component of $\partial M$, $\partial s$ is a covering map. 

If moreover $k \geq 1$, the  map $\partial s$  is conformal and its fibers have at most two elements.}
\end{enumerate}
\end{theorem}
So, under the hypothesis that $\partial_sM$ is not too bad (i.e it is locally the graph of a function) we have a very strong control on the topology of $\partial_sM$, and as soon as $\partial M$ is at least ${\cal C}^1$, the conformal structure induced by $h$ on $\partial_sM$ is also (almost completely) determined  by the conformal structure  of $\partial M$.
As an illustration, we get for example:

\begin{corollaire}
\label{hyperbolique}
Let $s$ be a smooth, strict conformal embedding of the $n$-dimensional hyperbolic space ${\bf H}^n$ into a smooth $n$-dimensional compact Riemannian manifold $(N,h)$, $n \geq 3$. If the boundary $\partial_s {\bf H}^n$, the topological boundary of the image $s({\bf H}^n)$, is locally the graph of a ${\cal C}^0$ function, it is a smooth hypersurface, which is conformally equivalent to the standard sphere ${\bf S}^n$ or the projective space ${\bf RP}^n$.
\end{corollaire}

This corollary follows directly from theorem \ref{geometrique}, because ${\bf H}^n$ is conformally diffeomorphic to the upper-hemisphere inside the round sphere ${\bf S}^{n}$.  So, inside ${\bf S}^{n}$, $\partial {\bf H}^n$ is a conformally flat $(n-1)$-sphere.  Now, if $s: {\bf H}^n \to (N,h)$ is a smooth strict conformal embedding, such that $\partial_s{\bf H}^n$ is locally the graph of a ${\cal C}^0$ function, theorem \ref{geometrique} says that $\partial_s{\bf H}^n$ is actually a smooth submanifold, which is conformally covered by the conformally flat sphere ${\bf S}^{n-1}$.  Since by the theorem, the cover is at most twofold, we get that  $\partial_s{\bf H}^n$ is either conformally diffeomorphic to ${\bf S}^{n-1}$, or to the $(n-1)$-dimensional projective space.

\subsection{Generalization to Cartan geometries}
\label{Cartan-intro}
Until there, we spoke only about Riemannian conformal structures. Nevertheless, a great part of the methods leading to the previous results are not at all specific to the Riemannian framework, but hold for more general geometric structures called {\it Cartan geometries}. Since this notion is central in all this work, and since it will allow us to get analogue results for very general geometric structures, it  seems  worthwhile to set the problem of geometric embeddings into the framewok of Cartan geometries.
 
Intuitively, a Cartan geometry is the data of a manifold infinitesimally modelled on some homogeneous space $\X=G/P$, where $G$ is a Lie group and $P$ a closed subgroup of $G$. Precisely, a Cartan geometry on a manifold $M$, modelled on the homogeneous space $\X=G/P$, is the data of:

- a principal $P$-bundle $\hm \rightarrow M$ over $M$.

- a $1$-form $\omega^M$ on $\hm$, with values in the Lie algebra $\lieg$, called {\it Cartan connexion}, and satisfying the following conditions:
\begin{enumerate}
\item{At every point $\hx \in \hm$, $\omega_{\hx}^M$ is an isomorphism between $T_{\hx}\hm$ and $\lieg$.}
\item{If $X^{\dagger}$ is a vector field of $\hm$, comming from the action by right multiplication of some one-parameter subgroup $t \mapsto Exp_G(tX)$ of $P$, then $\omega^M(X^{\dagger})=X$.}
\item{For every $a \in P$, ${R_a}^*\omega^M=Ad(a^{-1})\omega^M$ ($R_a$  standing for the right action of $a$ on $\hm$).}
\end{enumerate}
We will only consider Cartan geometries for which $P$ acts faithfully on $\X=G/P$, a condition which is satisfied in all cases of interest.

A lot of classical geometric structures, including pseudo-Riemannain metrics, pseudo-Riemannain conformal structures, $CR$, affine and projective structures  can be interpreted in terms of Cartan geometry, in the sense that they determine a Cartan geometry, which is canonical if one requires suitable normalization conditions on $\omega^M$ (see \cite{cap}, \cite{cartan}, \cite{chern}, \cite{kobayashi}, \cite{sharpe}, \cite{tanaka}). 

\subsubsection{Geometric embeddings}
\label{plongements-geometriques}
Let $\X=G/P$ be an homogeneous space, where $G$ is a Lie group and  $P$ a closed subgroup of $G$. We consider $(M,\hm,\omega^M)$ and  $(N,\hn,\omega^N)$ two Cartan  geometries  modelled on  $\X$ (notice that $M$ and $N$ automatically have the same dimension, which is that of $\X$). Such Cartan geometries will always supposed to be smooth. 

By a {\it geometric embedding} $\sigma$ of  $(M,\hm,\omega^M)$ into $(N,\hn,\omega^N)$, we mean a smooth bundle embedding $\sigma :\hm \to \hn$, such that  $\sigma^*(\omega^N)=\omega^M$. Such an embedding is equivariant for the right action of $P$ on $\hm$ and $\hn$ respectively, and induces a smooth  embedding $s: M \to N$.  We will say that  {\it $\sigma$ is a strict embedding} if $\sigma (\hm)$ is a strict open subset of $\hn$ (or equivalently $s(M)$ is a strict open subset of $N$). When the embedding $\sigma$ is not strict, we just have a  {\it geometrical isomorphism} between $(M,\hm,\om)$ and $(N,\hn,\on)$. In other words $(M,\hm,\om)$ and $(N,\hn,\on)$ are the same from the point of view of Cartan geometries modelled on $\X$.  

In the whole article, we will denote by  $\pshm$  (resp. $\psm$) the topological boundary of the open subset $\sigma(\hm)$ (resp. $s(M)$) in $\hn$ (resp. in $N$).

\subsubsection{Kleinian manifolds}
\label{presentation-kleiniennes}
Let $G$ be a Lie group, $P$ a closed subgroup of $G$, and $\X=G/P$. If $\omega^G$ denotes the Maurer-Cartan form on $G$,  the triple $(\X,G,\omega^G)$ is the {\it flat model} for the Cartan geometries modelled on $\X$. The general class of manifolds we will consider are {\it Kleinian manifold}, namely quotients $M=\Gamma \backslash \Omega$, where $\Omega$ is an open subset of $\X=G/P$, and $\Gamma$ a discrete subgroup of $G$, acting freely and properly on $\Omega$. As we will see in section \ref{kleinian-manifold}, a Kleinian manifold is naturally endowed with a canonical Cartan geometry modelled on $\X$, denoted $(M,\hm,\om)$. 

For Kleinian manifolds, we will prove an analogue of theorem \ref{codimension2}, namely:

\begin{theorem}
\label{Cartan-codimension}
Let $\Omega \subsetneq \X$ be an open subset such that  the Hausdorff dimension of the boundary $\partial \Omega$ is $< \text{dim}(\X)-1$. Let $M=\Gamma \backslash \Omega$ be a Kleinian manifold, and $(M,\hm,\om)$ its canonical Cartan geometry. Let $\sigma: (M,\hm,\om)  \to (N,\hn,\omega^N)$ be a geometric embedding, where  $(N,\hn,\omega^N)$ is a Cartan geometry modelled on $\X$. Then there is an open subset $\Omega^{\prime} \subset \XX$ containing $\Omega$, on which $\Gamma$ acts properly discontinuously, defining a Kleinian manifold $M^{\prime}={\Gamma\backslash\Omega^{\prime}}$, and such that $s$ extends to a geometric isomorphism $\sigma^{\prime} : (M^{\prime},{\hat M}^{\prime}, \omega^{M^{\prime}}) \to (N,\hn,\on)$.  \end{theorem}

As we said before, quite a lot of classical geometric structures are Cartan geometries, and among them, interesting examples are Kleinian manifolds. Thus, the previous theorem has nice illustrations, playing with different model spaces $\X$. To give a flavour of the kind of applications one can get, let us quote the:

\begin{corollaire}
\label{produit-hyp-sphere}
Let  $n \geq 3$, $n -1 > m \geq 1$ two integers, and ${\bf H}^m \times {\bf S}^{n-m}$ endowed with the conformal structure defined by the product of the hyperbolic and round metrics.  If $s : {\bf H}^m \times {\bf S}^{n-m} \to (N,h)$ is a conformal embedding into a Riemannian manifold of dimension $n$, then $(N,h)$ is conformally equivalent to an open subset of $({\bf S}^n,g_{can})$, whose boundary is an open subset of a conformal sphere of dimension $m-1$ in ${\bf S}^n$.
\end{corollaire}

The corollary follows directly from the theorem since ${\bf H}^m \times {\bf S}^{n-m}$ is conformally equivalent to the complement of a $(m-1)$-sphere in ${\bf S}^n$  (with the convention that a $0$-sphere consists of two points).

For open sets $\Omega$ bounded by a hypersurface, we will also prove:
\begin{theorem}
\label{klein-hypersurface}
Let $\Omega \subsetneq \X$ be an open subset such that $\overline \Omega$ is a topological submanifold with boundary of $\X$, and $\partial \Omega$ is of class ${\cal C}^{k,\alpha}$, $k+\alpha \geq 1$. Let $M=\Gamma \backslash \Omega$ be a Kleinian manifold, and $\sigma: (M,\hm,\om)  \to (N,B,\omega^N)$ be a geometric embedding, where  $(N,B,\omega^N)$ is a Cartan geometry modelled on $\X$. If the boundary $\partial_sM$ is locally the graph of a ${\cal C}^{0}$ map, then it is a hypersurface of $N$ of class ${\cal C}^{k,\alpha}$.
\end{theorem}

\subsection{Some results on maximal geometries}
\label{introduction-maximal}
The definition of geometric embeddings endows naturally the set of Cartan geometries modelled on a given $\X=G/P$ with a partial ordering.  If $(M,\hat M, \omega^M)$ and $(N,\hat N, \omega^N)$ are two such geometries, we will say that $M$ is smaller than $N$ if there is a strict embedding $\sigma : (M,\hat M, \omega^M) \to (N,\hat N, \omega^N)$. This leads to the following:

\begin{definition}(Maximal Cartan geometry)
\label{maximal-geometry}
A Cartan geometry $(M,\hm,\omega^M)$ modelled on  $\X$ is maximal when any geometrical embedding from  $(M,\hm,\omega^M)$ into another  Cartan geometry $(N,\hn,\omega^N)$ modelled on $\X$  is onto. 

\end{definition}

We already quoted results about conformally maximal Riemannian manifolds.  For Kleinian manifolds, we will prove:

\begin{theorem}
\label{klein-maximal}
Let $\Omega \subsetneq \X$ be a normal domain.  Let $M=\Gamma \backslash \Omega$ be a Kleinian manifold. If the action of $\Gamma$ is free and proper on no open subset of  $\overline \Omega$ containing strictly $\Omega$, then $M$ is maximal among the Cartan geometries modelled on $\X$.
\end{theorem}


The notion of normal domain will be introduced in section \ref{section-normal}.  Intuitively, a normal domain is an open set whose boundary is not too wild.  For example, open subsets whose boundary is a locally Lipschitz hypersurface, or have Hausdorff codimension $>1$ are normal.  The meaning of  theorem \ref{klein-maximal} is the following: if a Kleinian manifold is maximal among all Kleinian manifolds modelled on $\X$, it is maximal among (the larger class of) all Cartan geometries modelled on $\X$.

\subsection{Organization of the paper}
The article begins with a straigthforward generalization to general Cartan geometries of  the so called {\it $b$-boundary  construction}, first introduced by B.G Schmidt in \cite{schmidt2}. This construction associates a Cauchy boundary to any Cartan geometry. In section \ref{application-bord}, we explain how any strict geometrical embedding determines a boundary map from an open subset of this Cauchy boundary to the topological boundary of the image of the embedding.  In section \ref{quotient-manifold}, a more explicit description of the Cauchy boundary is provided for Cartan geometries obtained by quotienting open domains having a  nice enough topological boundary.  The  notion of weak Cartan geometry is then introduced. In section \ref{regularity}, an important regularity property of the boundary map is  proved (proposition \ref{extension}). Sections \ref{application-bord} and \ref{quotient-manifold} yield the general background  for the proofs of most of our results here, and we hope it will be also useful for further investigations on geometric embeddings (especially the case of pseudo-Riemannian (non Riemannian) conformal embeddings).  

Section \ref{riemannian-case} deals specifically with the case of conformal Riemannian structures, for which the general results obtained in sections \ref{bord} and \ref{quotient-manifold} can be sharpened.  This allows to prove theorems \ref{codimension2}, \ref{immersion} and \ref{geometrique}, using methods already introduced in \cite{frances1}.  This section is the most involved part of the paper.

In section \ref{kleinian-manifold}, we investigate geometrical embeddings of Kleinian manifolds, proving theorems \ref{Cartan-codimension} and \ref{klein-hypersurface}.  Section \ref{section-maximal} is devoted to conformally maximal manifolds, and more generally maximal Cartan geometries.  The reader will find here the proofs of theorem \ref{klein-maximal}, from which theorem \ref{courbure-constante-maximale} is deduced.  Illustrations in the framework of $CR$-structures, affine, Lorentz  and projective geometries are given in \ref{applications}. We also prove in  section \ref{section-maximal} theorems \ref{homogene} and \ref{gros-isometries}, using tools introduced in section \ref{riemannian-case}.

\section{ The boundary map of a geometric embedding}
\label{bord}
\subsection{Cauchy completion of a Cartan geometry}
\label{bord-Cauchy}
We consider a Cartan geometry $(M,\hm,\om)$ modelled on the space $\X=G/P$.
Let us fix once for all  $X_1,...,X_m$ a basis of the  Lie algebra $\lieg$. This choice defines a parallelism  ${\cal R}$ on $\hm$ as follows: 
$${\cal R}_{\hx}=((\omega_{\hx}^M)^{-1}(X_1),...,(\omega_{\hx}^M)^{-1}(X_m)), \ \ \hx \in \hm.$$

 There is a unique Riemannian metric  $\rho^M$ on  $\hm$ for which this parallelism is orthonormal, and we denote by  $d_M$  the  distance defined by $\rho^M$ on $\hm $. 
 Let $({\hat M}_c, {\overline d}_M)$ denote the Cauchy completion of the metric space  $(\hm,d_M)$, and  $\partial_c\hm$ its  Cauchy boundary, i.e  $\partial_c\hm = \hm_c \setminus \hm$. Let us remark that for every $p \in P$, the Jacobian matrix of $R_p$, the right multiplication by $p$, when  expressed in the parallelism ${\cal R}$, is exactly the matrix of the linear map $\Ad p^{-1}$ expressed in the basis $(X_1,...,X_n)$. It follows that the maps $R_p$ are uniformly continuous, thus extend to continuous  maps ${\overline R}_p$ on $\hm_c$, which  leave  $\partial_c\hm$ invariant. As a consequence, one can also define the {\it Cauchy completion} ${M_c}$ (resp. the {\it Cauchy boundary}  $\partial_cM$ ) of $M$ as the quotient $\hm_c/P$ (resp. $\partial_c\hm/P$). 
 
Observe that the initial choice of the basis $(X_1,...,X_m)$ is not relevant, since  different choices lead to bi-Lipschitz equivalent metrics, hence to the same Cauchy completion.

The space ${M_c}$ is generally quite far from being  a nice topological space, because the action of $P$, which is proper on $\hm$, can behave very badly on the boundary $\partial_c\hm$, yielding a non Hausdorff quotient.
\subsubsection{Two easy examples}
\label{illustration-kleinienne}

As a first example, let us determine the Cauchy boundary of the model space $\X$ itself.  The Cartan geometry that we are considering here is the triple $(\X,G,\omega^G)$, where $\omega^G$ is the Maurer-Cartan form on $G$.  Since the Maurer-Cartan form is  invariant by left multiplication, the left action of $G$ leaves the metric $\rho^{\X}$ invariant, so that $(G,\rho^{\X})$ is an homogeneous Riemannian manifold, hence complete. As a consequence $\partial_c{\hat \X}=\partial_cG=\emptyset=\partial_c \X$.

Another easy, and maybe more enlightening example is that of $M=\re^n$, $n \geq 3$, endowed with its standard conformally flat structure.  Here, the model space  $\X$ is the sphere ${\bf S}^n$  with its  standard conformal structure, and seen as the homogeneous space  $G/P$, where $G=SO(1,n+1)$ and $P$ is a parabolic subgroup. Since   the stereographic projection identifies conformally $\re^n$ with the sphere minus a point, the bundle  $\hm$ is identified with the open subset of $SO(1,n+1)$, obtained  by removing a $P$-orbit for the right action. The Cauchy boundary $\pcho$  is identified with the topological boundary $\pho$, i.e a  $P$-orbit for the right action. We thus get that  $\partial_c\re^n$ is just a point. 

For the conformal structure defined by a flat complete Riemannian  manifold $M=\Gamma \backslash \re^n$, with $\Gamma \subset Is(\re^n)$ nontrivial and discrete, one checks that $\partial_c\hm$ is identified with $\Gamma \backslash P$, so that $\partial_cM$ is still a point. Nevertheless, when $\Gamma$ is noncompact, the space $M_c$ is not Hausdorff anymore.





\subsection{Construction of the boundary map}
\label{application-bord}
Let $(M,\hm,\omega^M)$ and $(N,\hn,\omega^N)$ be two Cartan geometries modelled on the same homogeneous space $\X=G/P$, and $\sigma : (M,\hm,\om) \to (N,\hn,\on)$ be a geometric embedding. In what follows, we will call $s$ the embedding from $M$ into $N$ induced by $\sigma$.  By, $\partial_{\sigma}\hm$  (resp.  $\partial_sM$), we will understand the topological boundary of $\sigma(\hm)$ in $\hn$ (resp.  of $s(M)$ in $N$). The construction described in section \ref{bord-Cauchy}, {\it using a same fixed basis ($X_1,...,X_m$) of $\lieg$}  yields two Riemannian metrics $\rho^M$ and $\rho^N$ on $\hm$ and $\hn$ respectively. The construction implies $\sigma^*(\rho^N)=\rho^M$,  namely $\sigma$ is {\it an isometric embedding} from $(\hm,\rho^M)$ into $(\hn,\rho^N)$. Let  us define $d_{N}^{\sigma}$ as  the length distance induced by  $\rho^N$ on $\sigma(\hm)$. If $\hy$ and  $\hz$ are in $\sigma(\hm)$, then:
 $$d_{N}^{\sigma}(\hy,\hz)=\inf\{  L_{\rho^N}(\gamma) \ | \ \gamma : [0,1] \to \sigma(\hm) \ of \ class \ C^1, \ \gamma(0)=\hy, \ \gamma(1)=\hz \}$$

 We always have $d_{N}^{\sigma} \geq d_N$, what will be often used in the following to say that a  Cauchy sequence for  $d_{N}^{\sigma}$ is automatically a  Cauchy sequence for $d_N$. The map $\sigma$ is an isometry between $(\hm,d_M)$ and $(\sigma(\hm),d_{N}^{\sigma})$.




\begin{definition}
\label{cauchy-regular-set}
Let $\sigma : (M,\hm,\om) \to (N,\hn,\on)$ be a geometric embedding. The regular set of $\partial_c\hm$ associated to $\sigma$ is defined as:
 $$\hl_c = \{ \hx \in \pchm \ |  \ \exists  (\hx_k) \ a \ sequence \ of \  \hm,  \ \hx_k \to \hx \ and \  \sigma(\hx_k) \ converges \ in \ \hn \}$$ 
\end{definition}

Since $\sigma$ is equivariant for the action of $P$ on $\hm$ and $\hn$, the set $\hl_c$ is $P$-invariant and we call $\Lambda_c=\hl_c/P$. The importance of the set $\hl_c$ comes from the fact that there is a nice notion of boundary map defined on it.  Before making it precise, let us introduce the:

\begin{definition}
\label{Accessible points}
A point $x \in \psm$ is said to be accessible if there exists a $C^1$ path  $\gamma : [0,1] \to N$ such that $\gamma([0,1[) \subset s(M)$ and  $\gamma(1)=x$. A point in  $\pshm$ is called accessible if it is in the fiber of an accessible point of $\psm$. 
\end{definition}

It is clear that  $\hx \in \pshm$ is accessible iff there exists a $C^1$ path   $\hat \gamma : [0,1] \to \hn$ such that  $\hat \gamma ([0,1[) \subset \sigma(\hm)$ and $ \hat \gamma (1)=\hx$.

We can now state: 

\begin{proposition}
\label{boundary}
Let us assume that the embedding $\sigma$ is strict. Then:
\begin{enumerate}
\item{The set $\hl_c$ is a nonempty open subset of $\partial_c\hm$.}
\item{If $\hx \in \hl_c$, then for any sequence $(\hx_k)$ of $\hm$ tending to $\hx$, the sequence $\sigma(\hx_k)$ converges. The limit, denoted $\partial \sigma(\hx)$, depends only of $\hx$, what yields a well defined map $\partial \sigma : \hl_c \to \pshm$.}
\item{ The map $\overline \sigma: \hl_c \cup \hm \to \pshm \cup \sigma(\hm)$ which coincides with $\sigma$ on $\hm$, and with  $\partial \sigma$ on $\pchm$ is  $P$-equivariant and continuous. The maps $\partial \sigma$ and $\overline \sigma$ induce continuous maps $\partial s: \Lambda_c \to \partial_s M$ (called the boundary map of $s$) and $\overline s : \Lambda_c \cup M \to \partial_sM \cup s(M)$. } 
\item{The image $\partial s (\Lambda_c)$  (resp.  $\partial \sigma (\hl_c)$) contains every accessible point of $\partial_sM$  (resp.  of $\pshm$).  In particular, it is dense in $\partial_sM$  (resp.  in $\pshm$).}
\item{The  group $P$ acts freely and properly on  $\hl_c \cup \hm$; in particular  $\Lambda_c \cup M \subset M_c$ is Hausdorff.}
\end{enumerate}

\end{proposition}

\begin{preuve}
We use the following lemma, which  is probably standard (see \cite{frances1}, Lemma $14$):

\begin{lemme}
\label{accessible}
When  $\psm$ is nonempty, the set of accessible points is dense in $\psm$.
\end{lemme}

Since we assumed that $s$ is strict, $\psm$ is nonempty, and by the previous lemma, the accessible points are dense in $\psm$.  Let us choose  $\hz$ accessible in  $\pshm$. By definition, there is a $C^1$ path $\gamma : [0,1] \to \hn$ such that $\gamma([0,1[ ) \subset \sigma(\hm)$, and  $\gamma(1)=\hz$.  Let $(t_k)$  be a sequence of  $[0,1[$ tending to $1$. Since the distance $d_N^{\sigma}(\gamma(t_k), \gamma(t_{k+p}))$ is always less than the $\rho^N$-length of the segment $\gamma([t_k,t_{k+p}])$, we get that $\gamma(t_k)$ is a Cauchy sequence for  $d_N^{\sigma}$. Let $\hx_k$ be the point of $\hm$ such that  $\sigma(\hx_k)=\gamma(t_k)$. Then $(\hx_k)$ is a  Cauchy sequence for $d_M$, hence converges to $\hx_{\infty} \in \pchm$. By construction $\hx_{\infty} \in \hl_c$, what proves that  $\hl_c$ is nonempty. Let now $\epsilon$ be so small that $\overline B^N(\hz,\epsilon)$ (the closed $d_N$-ball of center $\hz$ and radius $\epsilon$) is complete for the distance $d_N$.  Let $\hy \in B^M(\hx_{\infty},\frac{\epsilon}{2})$ (the open ${\overline d}_M$-ball of center  $\hx_{\infty}$ and radius $\frac{\epsilon}{2}$), and  $(\hy_k)$ a sequence of  $\hm$ tending to $\hy$. Then $\sigma(\hy_k)$ is a Cauchy sequence for $d_N^{\sigma}$, hence for $d_N$, and for $k$ big enough: $\sigma(\hy_k) \in \overline B^N(\hz,\epsilon)$.  We get that  $\sigma(\hy_k)$ converges in $\hn$,    and  $\hy \in \hl_c$. As a consequence, $\hl_c$ is open in $\pchm$.

Let us now prove the second point of the proposition.  Let  $\hx$ be a point of $ \hl_c$. By definition, there is a sequence $(\hx_k)$ of $\hm$ converging to $\hx$, such that  $\sigma(\hx_k)$ converges to some point of  $\partial_{\sigma}\hm$.  If $(\hx_k^{\prime})$ is another sequence of $\hm$ tending to  $\hx$, then $d_M(\hx_k,\hx_k^{\prime})=d_{N}^{\sigma}(\sigma(\hx_k),\sigma(\hx_k^{\prime}))$ tends to $0$. {\it A fortiori} $d_N(\sigma(\hx_k),\sigma(\hx_k^{\prime}))$ tends to $0$, what means that $\sigma(\hx_k^{\prime})$ converges to the same limit as $\sigma(\hx_k)$. This shows that  ${\partial \sigma}(\hx)$ is well defined. It is clear that $\overline \sigma$ is $P$-equivariant, and it is continuous since it is $1$-Lipschitz with respect to the distances ${\overline d}_M$ and $d_N$.  The third point of the proposition comes from the $P$-equivariance of $\partial \sigma$ and $\overline \sigma$.

At the begining of the proof, we showed that if $\hz$ is accessible in $\pshm$, there is a sequence  $(\hx_k)$ of $\hm$ tending to $\hx \in \hl_c$, such that $\sigma(\hx_k)$ converges to $\hz$.  This exactly means  $\partial \sigma(\hx)=\hz$, what proves that $\partial \sigma(\hl_c)$ contains every accessible point of $\pshm$.  By $P$-equivariance, the same property is true for $\partial s$, and this gives the fourth point of the proposition.

To prove the last point of the proposition,  consider $(\hx_k)$  a sequence of  $\hl_c \cup \hm$ converging to $\hx_{\infty} \in \hl_c \cup \hm$, and $(p_k)$ a sequence of $P$ such that $\hx_k.p_k$ converges  to $\hy_{\infty} \in \hl_c \cup \hm$.   Then  $\overline \sigma(\hx_{k})$ tends to $\overline \sigma(\hx_{\infty})$ and $\overline \sigma(\hx_{k}).p_k$ tends to $\overline \sigma(\hy_{\infty})$. Since $P$ acts properly on  $\hn$, the sequence  $(p_k)$ is bounded in $P$. This shows that the action of  $P$ is proper on $\hlc \cup \hm$. The action of  $P$ is also free on  $\hlc \cup \hm$, because it is free on  $\hn$, and if  $\hx.p=\hx$, for some  $\hx \in \hlc \cup \hm$ and $p \in P$, then $\overline \sigma (\hx).p=\overline \sigma (\hx)$.

\end{preuve}

\subsection{First consequences of proposition \ref{boundary}}

It is worth noticing that proposition \ref{boundary} has interesting implications.  For example, if $(M,\hm,\om)$ is a Cartan geometry such that $\partial_cM$ has only finitely many points, and if $\sigma :  (M,\hm,\om) \to (N,\hn,\on)$ is a strict geometrical embedding, inducing an embedding $s: M \to N$.   Then $\Lambda_c$ also has finitely many points, while $\partial s(\Lambda_c)$ must be  dense in $\partial_sM$.  This implies that $\partial_sM$ must have finitely many points.  An extreme case is that of the conformal structure of the Euclidean space, for which we saw that $\partial_c \re^n$ is only one point.  Then, we get that if $s : (\re^n,g_{eucl}) \to (N,h)$ is a strict conformal embedding into some $n$-dimensional Riemannian manifold $(N,h)$, $n \geq 3$, then $\partial_s\re^n$ must be a single point, and then, it is easy to see that $(N,h)$ is conformally equivalent to the round sphere.

Another  straigthforward consequence of proposition \ref{boundary}, which will be exploited in the proof of theorem \ref{klein-maximal}, is that whenever $P$ does not act properly and freely on any nonempty open subset of $\partial_cM$, the Cartan geometry $(M,\hm,\om)$ has to be  geometrically maximal.






\section{Normal domains, quotient manifolds and  weak Cartan geometries}
\subsection{Quotient manifolds}
\label{quotient-manifold}
We consider $(L,\hL, \omega^L)$ a Cartan geometry modelled on $\X=G/P$.  If $\Omega$ is an open subset of $L$, then $\Omega$ is naturally endowed with a Cartan geometry modelled on $\X$, induced by that of $L$.  The Cartan bundle $\hat \Omega$ is the inverse image of $\Omega$ under the projection $\pi_L: \hL \to L$. The Cartan connexion $\oo$ on $\ho$ is just the restriction of $\ol$ to $\ho$.
As we saw in section \ref{bord-Cauchy}, the choice of a basis $X_1,...,X_m$ of $\lieg$ endows $\hL$ (resp. $\ho$) with a Riemannian metric $\rho^L$ (resp. $\rho^{\Omega}$). The metric $\rho^{\Omega}$ is just the restriction of $\rho^L$ to $\ho$.  We call $d_L$ (resp.  $d_{\Omega}$) the distance defined by $\rho^L$ on $\hL$ (resp. defined by $\rho^{\Omega}$ on $\ho$).

We say that a diffeomorphism $\phi$ of $L$ is an automorphism of the Cartan geometry if it lifts to $\hat \phi$, a bundle automorphism of $\hat L$ satisfying ${\hat \phi}^* \omega^L=\omega^L$. The group of automorphisms of $L$ will be denoted by $\Aut L$. The condition that $P$ acts faithfully on $\X=G/P$ implies that an element of $\Aut L$ lifts in an unique way to a bundle automorphism of $\hat L$.
  Assume now that $\Gamma$ is a discrete subgroup of $\Aut L$, preserving $\Omega$ and acting freely and properly on it. Let us call $M$ the quotient manifold $\Gamma \backslash \Omega$. Such a manifold also inherits from $L$ a Cartan geometry modelled on $\X$, that we decribe now.  
  
\subsection{Weak Cartan geometries and canonical quotient geometries}
\label{weak}
Because it preserves a parallelism, $\Aut L$ acts freely and properly on $\hL$, and its action commutes with that of $P$. We denote by $\hlg$ the quotient manifold $\Gamma \backslash \hL$, and by $\omega^{\Lg}$ the $1$-form induced by $\omega^{L}$ on $\hlg$. The manifold $\hlg$ still carries a right action of $P$, and the $1$-form $\omega^{\Lg}$ satisfies the conditions $1$,$2$ and  $3$ of a Cartan geometry (see section \ref{Cartan-intro}).  Nevertheless, $\hlg$ is generally not a $P$-principal bundle (for example, unless the action of $\Gamma$ is proper on the whole $L$, the action of $P$ will not be proper on $\hlg$). Thus,  $(\hlg,\omega^{\Lg})$ is not a Cartan geometry in the classical sense:  we will call it  {\it a weak Cartan geometry}. 

If $\hm$ is the projection of $\ho$ on $\hlg$, then $\hm$ is stable for the right $P$-action on $\hlg$, and is a $P$-principal bundle over $M$.  If $\omega^M$ denotes the restriction of $\omega^{\Lg}$ to $\hm$, then the triple $(M,\hm,\omega^M)$ is a Cartan geometry.  We call it {\it the canonical Cartan geometry on $M$, induced by that of $L$}.

\subsection{Normal domains}
\label{section-normal}
The metric $\rho^{L}$ (resp.  $\rho^{\Omega}$) is $\Gamma$-invariant, hence induces a metric $\rho^{\Lg}$ (resp. $\rho^M$) on $\hlg$ (resp. on $\hm$), as well as a distance $d_{\Lg}$ (resp. $d_M$). Observe that the metric $\rho^M$ is precisely the one constructed in section \ref{bord-Cauchy} for the Cartan geometry $(M,\hm,\om)$.  Observe also that the projection $\pi_{\Gamma} : (\hL, \rho^L)  \to (\hlg, \rho^{\Lg})$ is a Riemannian covering. 
The previous construction leads naturally to the following problem:

{\it What is the link between the Cauchy boundary  $\partial_c\hm$ of the metric space $(\hat M,d_M)$, as described in section \ref{bord-Cauchy}, and the  topological boundary $\partial \hm$ in $\hlg$?}

When $\partial \hm$ is very irregular, the relationship between the two boundaries is not clear, because a sequence of $\hm$ may converge to a  point of $\partial \hm$ without being a Cauchy sequence for $d_M$.  Nevertheless, we are going to exhibit a class of open sets for which the link between  $\partial_c\hm$ and $\partial \hm$ is quite easy to understand.

\begin{definition}[Normal domains]
\label{normal-domains}
Let $Y$ be a connected smooth  manifold and $W \subsetneq Y$ an open subset.  We say that $W$ is a normal domain of $Y$ if for any $y \in \partial W$, there exists a countable family of connected relatively compact neighbourhoods ${\cal U}_y=\{  U_i\}_{i \in \NN}$  such that:
\begin{enumerate}
\item{For every $i \in \NN$, ${ \overline U}_{i+1} \subsetneq U_i$, and $\bigcap_{i \in \NN} U_i=\{ y \}$}
\item{For any neighbourhood $U \in {\cal U}_y$, $U \cap W$ is connected}
\item{For every smooth Riemannian metric $\rho$ on $U_0$, and every $i >1$, the metrics $d_{\rho}^{U_i}$ and $d_{\rho}^{U_i  \cap W}$ are bi-Lipschitz equivalent on $U_i \cap W$.  Namely, there exists a real $k_i \geq 1$ such that for any $y_1,y_2 \in U_i \cap W$:
$$ d_{\rho}^{U_i}(y_1,y_2) \leq d_{\rho}^{U_i \cap W}(y_1,y_2) \leq k_i d_{\rho}^{U_i}(y_1,y_2)$$}

\end{enumerate}

\end{definition}
Recall that $d_{\rho}^{U_i}$ and $d_{\rho}^{U_i \cap W}$ are defined as follows:
$$ d_{\rho}^{U_i }(y_1,y_2)=\inf \{ L_{\rho}(\gamma) \ | \  \gamma \in C^1([0,1] ; U_i ), \  \gamma(0)=y_1 \  {\sc and}  \ \gamma(1)=y_2 \}$$
$$ d_{\rho}^{U_i \cap W }(y_1,y_2)=\inf \{ L_{\rho}(\gamma) \ | \  \gamma \in C^1([0,1] ; U_i \cap W), \  \gamma(0)=y_1 \  {\sc and}  \ \gamma(1)=y_2 \}$$

\begin{remarque}
\label{remarque-normal}
The property that $W$ is a normal domain does not depend on the choice of the smooth Riemannian metric $\rho$ on $U_0$, because two such Riemannian metrics define bi-Lipschitz equivalent distances on each $U_i$ (resp.  each $U_i \cap W$), $i \geq 1$.  So, to show that  $W \subset Y$ is a normal domain, it is sufficient to check that it satisfies the previous definition for one  smooth Riemannian metric on $U_0$.
\end{remarque}

The following lemma exhibits a wide class of normal domains.  Its proof, as well as that of lemma \ref{normal-fibre}, is postponed at the end of the paper (section \ref{annexe}).

\begin{lemme}
\label{exemples-normaux}
Let $Y$  be a smooth manifold, and $W \subset Y$   a strict open subset. Assume that we are in one of the following cases:
\begin{enumerate}
\item{$\overline W$ is a topological submanifold with boundary of $Y$, and $\partial W$ is locally the graph of a Lipschitz function.}
\item{$\partial W$ has Hausdorff dimension $< n-1$.}
\end{enumerate}
Then $W$ is a normal domain.

\end{lemme}

In the following, we will also need:

\begin{lemme}
\label{normal-fibre}
Let $Y$ be a smooth manifold, and $\pi : B \to Y$ a smooth fiber bundle. Assume that  $W \subset Y$ is a normal domain, and define $\check W \subset B$ as $\check W : = \pi^{-1}(W)$.  Then $\check W$ is a normal domain of $B$.

\end{lemme}


\subsection{The natural map $\lambda : \partial \hm \to \partial_c\hm$}
\label{natural-map}
Let $(L,\hL,\ol)$ be a Cartan geometry modelled on $\X$, of dimension $n$. Let  $\Omega \subset L$ be  a strict open subset and $M=\Gamma \backslash \Omega$ a quotient manifold, where $\Gamma$ is discrete in $\Aut L$.  We assume now that $\Omega \subsetneq L$ is a {\it normal domain}. By lemma \ref{normal-fibre}, $\hat \Omega$ is a normal domain of $\hL$.  The projection from $\hat L$ to $\hL_{\Gamma}$ is a local isometry (with respect to the metrics $\rho^L$ and $\rho^{L_{\Gamma}}$), so that  if $M=\Gamma \backslash \Omega$, then $\hm$ is a normal domain of $\hL_{\Gamma}$. In particular, to each $\hx \in \phm$  is associated a real $1 \leq k_{\hx}$, and a  neighbourhood ${ U}_{\hx}$ such that $U_{\hx} \cap \hm$ is connected, and satisfying: 
\begin{equation}
\label{equivalence}
 d_{L_{\Gamma}}^{U_{\hx}}(\hx_1,\hx_2) \leq d_{L_{\Gamma}}^{U_{\hx} \cap \hm}(\hx_1,\hx_2) \leq k_{\hx} d_{L_{\Gamma}}^{U_{\hx}}(\hx_1,\hx_2)
 \end{equation}

Recall that  $d_{L_{\Gamma}}^{U_{\hx} \cap \hm}$ is the distance induced by $\rho^{\Lg}$  on $U_{\hx} \cap \hm$, namely:
$$ d_{L_{\Gamma}}^{U_{\hx} \cap \hm}(\hx_1,\hx_2)=\inf \{ L_{\rho^{\Lg}}(\gamma) \ | \  \gamma \in C^1([0,1] ; U_{\hx} \cap \hm ), \  \gamma(0)=\hx_1 \  {\sc and}  \ \gamma(1)=\hx_2 \}.$$

Let $\hx$ be a point of $\phm$, and $(\hx_k)$ a sequence of $U_{\hx} \cap \hm$ converging to $\hx$.  The inequality (\ref{equivalence}) shows that  $(\hx_k)$, which is a Cauchy sequence for $d_{\Lg}^{U_{\hx}}$, is also a Cauchy sequence for $d_{L_{\Gamma}}^{U_{\hx} \cap \hm}$, hence for $d_M$, because  $d_{M} \leq d_{L_{\Gamma}}^{U_{\hx} \cap \hm}$. Thus the data of $\hx$ and $(\hx_k)$ naturally defines a point of $\partial_c\hm$. Now, if $(\hx_k^{\prime})$ is another sequence of $U_{\hx} \cap \hm$ converging to $\hx$, then it is equivalent to $(\hx_k)$ for $d_{\Lg}^{U_{\hx}}$,  hence for $d_{L_{\Gamma}}^{U_{\hx} \cap \hm}$ by  (\ref{equivalence}), hence for $d_M$.  Thus, the single data of a point in $\partial \hm$ defines a point of $\partial_c\hm$. This yields a well defined $P$-equivariant map $\lambda : \phm \to \partial_c \hm$.

\begin{lemme}
\label{bord-normal}
Let $(L,\hL,\ol)$, $\Omega$ and $M=\Gamma \backslash \Omega$ be as above.  Then the map $\lambda : \phm \to \partial_c\hm$ maps $\phm$ homeomorphically onto an open subset of  $\partial_c\hm$

\end{lemme}

\begin{preuve}
Let $\hx \in \phm$, and ${U}_{\hx}$ as above.   On $(U_{\hx} \cap \phm,d_{L_{\Gamma}}^{U_{\hx}})$, the map $\lambda$ is $k_{\hx}$-Lipschitz, hence continuous. It is injective because $d_{L_{\Gamma}} \leq d_M$. To show that the image of $\lambda$ is open, let us choose $\epsilon$  small enough so that the closed $\rho^{\Lg}$-ball of center $\hx$ and radius $\epsilon$ is complete for $d_{L_{\Gamma}}$, and included in $U_{\hx}$. Let $(\hx_k)$ be, as above,  a Cauchy sequence of $\hm$ (for the distance $d_M$) converging to $\hx$.  Now, let $\hy^{\prime} \in \partial_c\hm$ such that ${\overline d}_M(\lambda(\hx),\hy^{\prime}) < \frac{\epsilon}{2}$. If $(\hy_k^{\prime})$ is a  Cauchy sequence  of $(\hm,d_M)$ tending to $\hy^{\prime}$, then   ${d}_{M}(\hy_k^{\prime},\hx_k) \leq \frac{ \epsilon}{2}$ for $k$ big enough.  This implies that ${d}_{L_{\Gamma}}(\hy_k^{\prime},\hx_k) \leq \frac{ \epsilon}{2}$, so that $(\hy_k^{\prime})$ is a Cauchy sequence for $d_{\Lg}$, which is, for $k$ big enough, included in the the closed $\rho^{\Lg}$-ball of center $\hx$ and radius $\epsilon$. Thus, it must converge to $\hx^{\prime } \in \partial \hm \cap U_{\hx}$, and we get $\hy^{\prime} = \lambda(\hx^{\prime}) $. 
 
Starting with ${\overline d}_M(\lambda(\hx),\lambda(\hx^{\prime})) < \frac{\epsilon}{2}$, we obtained  $d_{L_{\Gamma}}(\hx,\hx^{\prime}) \leq \epsilon$. Hence, the map $\lambda^{-1}$ is continuous, and $\lambda$ is an  homeomorphism on its image.
\end{preuve}


Now, assume that $\Omega \subset L$ is a normal domain, and that $M=\Gamma \backslash \Omega$,  with $\Gamma$ discrete in $\Aut L$. 
Assume that $(N,\hn,\on)$ is a Cartan geometry modelled on $\X$, and that $\sigma : (M,\hm,\om) \to (N,\hn,\on)$ is a strict geometrical embedding (here, $(M,\hm,\om)$ is the natural Cartan geometry induced by $L$ on $M$, as described previously). Thanks to the map $\lambda$, we identify $\partial \hm$ with an open subset of $\partial_c\hm$, and we would like to understand  the subset $\hl=\hl_c \cap \partial \hm$, namely:
$$\hl= \{  \hx \in \phm \ | \  \exists (\hx_k) \ a \ sequence \ of \ \hm, \ \hx_k \to \hx \ and \ \sigma(\hx_k) \ converges \ in \ \hn \}$$
This will be done thanks to the:

\begin{lemme}
\label{identification}
If $\hl \not = \emptyset$, then $(\hl \cup \hm)/P$  is a Hausdorff open subset of $\Lambda_c \cup M$, homeomorphic to $\Gamma \backslash (\tilde \Lambda \cup \Omega)$, for $\tilde \Lambda \cup \Omega$ an open subset of $\overline \Omega$, containing $\Omega$ strictly, and  on which $\Gamma $ acts freely and properly.
\end{lemme}

\begin{preuve}
We already saw in  Proposition \ref{boundary} that $P$ acts freely and properly on $\hl_c \cup \hm$. Since $\hl \cup \hm$ is $P$-equivariantly homeomorphic to an open subset of $\hl_c \cup \hm$,   $\Lambda \cup M= (\hl \cup \hm)/P$ is Hausdorff. Let us call $\check  \Lambda \cup \hat \Omega$ the inverse image of $\hl \cup \hm$ by the projection $\hL \to {\hat L}_{\Gamma}$, and $\tilde \Lambda \cup \Omega$ the projection of $\check \Lambda \cup \hat \Omega$ on $\overline \Omega$. Then it is not difficult to check that since $P$ acts freely properly on $\hl \cup \hm$, $\Gamma$ acts freely properly on $\tilde \Lambda \cup \Omega$, and that $\Lambda \cup M$ is the quotient $\Gamma \backslash (\tilde \Lambda \cup \Omega)$.   
\end{preuve}




\subsection{Regularity properties for the boundary map}
\label{regularity}

Since $\hl$ and $\Lambda=\hl/P$ are identified with open subsets of $\hlc$ and $\Lambda_c$ respectively, we can restrict the boundary maps $\partial \sigma$ and $\partial s$
to them. If $\hx \in \hl$, (resp. $x \in \Lambda$) then $\partial \sigma(\hx)$ (resp.  $\partial s(x)$) is the limit in $\partial_{\sigma}\hm$  (resp. in $\partial_sM$) of $\sigma(\hx_k)$ (resp.  $s(x_k)$), for any sequence $(\hx_k)$ (resp.  $(x_k)$) of $\hm$ (resp. of $M$) tending to $\hx$ (resp. to $x$).  

We now state a proposition, which says that around every point of $\hl$, the map $\overline \sigma$ is the restriction of a smooth diffeomorphism. This will be crucial in the following, since $\partial \sigma$ will thus enjoy nice regularity conditions as soon as $\partial \hm$ does.
\begin{proposition}
\label{extension}
 Assume that $L$, $\Omega$, $M$ and $\Gamma$ satisfy the hypotheses of Lemma \ref{exemples-normaux}. Then for every $\hx \in \hl$, there is a neighbourhood $V_{\hx}$ (resp $V_{\hy}$) of $\hx$ in $\hlg$ (resp. of $\hy=\partial \sigma(\hx)$ in $\hn$), and $\sigma_{\hx}^{\prime}$ a smooth diffeomorphism between $V_{\hx}$ and $V_{\hy}$, such that $\sigma_{\hx}^{\prime}$ restricts to $\overline \sigma$ on $V_{\hx} \cap (\hl \cup \hm)$.
\end{proposition}

Of course, the proposition above is useless if $\hl= \emptyset$. Nevertheless, we will exhibit quite a lot of interesting examples (for instance the case of conformal Riemannian embeddings), where we will be able to show that $\hl \not = \emptyset$.

\begin{preuve}
We first make the assumption that the Hausdorff dimension of $\partial \Omega$ is $< n-1$. For any $\hx \in \hlg$ (resp $\hy \in \hn$ ), we will denote by $Exp_{\hx}$ (resp. $Exp_{\hy}$) the exponential map defined by the Riemannian metric $\rho^{\Lg}$ (resp. $\rho^{N}$). Let $\hx \in \hl$, and $(\hx_k)$ a sequence of $\hm$ tending to $\hx$, such that $\hy_k=\sigma(\hx_k)$ tends to $\hy \in \hn$.  For $k_0$ sufficiently large, there is $U_{\hx_{k_0}} \subset T_{\hx_{k_0}}\hlg$ a neighbourhood of $0_{\hx_{k_0}}$, such that $Exp_{\hx_{k_0}}(U_{\hx_{k_0}})=V_{\hx}$ is a simple convex neighbourhood for the metric $\rho^{\Lg}$, which  contains $\hx$. This means that any two points of $V_{\hx}$ can be joined by a unique geodesic segment (for the metric $\rho^{\Lg}$) included in $V_{\hx}$. We will also assume that we have choosen $k_0$ big enough, and $U_{\hx_{k_0}}$ small enough, so that $Exp_{\hy_{k_0}}$ is a diffeomorphism from $U_{\hy_{k_0}}=(\omega_{\hy_{k_0}}^N)^{-1} \circ \omega_{\hx_{k_0}}^{\Lg}(U_{\hx_{k_0}})$ on its image $V_{\hy}$.  Now, define $\sigma_{\hx}^{\prime} : V_{\hx} \to V_{\hy}$ by: 
$$\sigma_{\hx}^{\prime}=Exp_{\hy_{k_0}}\circ D_{\hx_{k_0}}\sigma \circ Exp_{\hx_{k_0}}^{-1}$$

The map $\sigma_{\hx}^{\prime}$ is a smooth diffeomorphism between $V_{\hx}$ and $V_{\hy}$.  Now, because $\sigma$ is an isometry between $(\hm,\rho^M)$ and $(\sigma(\hm),\rho^N)$, we have that $\sigma_{\hx}^{\prime}$ coincides with $\sigma$ on $W_{\hx_{k_0}}$, the set of points in $V_{\hx} \cap \hm$ which can be joigned to $\hx_{k_0}$ by a geodesic segment included in $V_{\hx} \cap \hm$. Now, the codimension asumption on $\partial \Omega$ allows  to prove the:

\begin{lemme}
\label{dense}
For every $\hz \in V \cap \hm$, the set $W_{\hz}$ of points in $V \cap \hm$, which can be joigned to $\hz$ by a $\rho^{\Lg}$-geodesic segment included in $V \cap \hm$ is dense and open in $V$.

\end{lemme}

\begin{preuve}
Since $V$ is supposed to be a convex, simple neighbourhood, for every $\hz \in V$, there is $V_{\hz} \subset T_{\hz}\hlg$, such that $Exp_{\hz}$ is a diffeomorphism from $V_{\hz}$ onto $V$.  Observe that $V_{\hz}$ is star-shaped relatively to $0_{\hz}$.  Now, let $S_{\hz}(r)$ (resp.  $B_{hz}(r)$) be the sphere (resp. the open ball) of radius $r$ and center $0_{\hz}$, for the metric $\rho_{\hz}^{\Lg}$. When $r$ is small enough, $S_{\hz}(r) \subset V_{\hz}$, and there is a well defined radial smooth projection $\pi : V_{\hz} \setminus \{ 0_{\hz} \} \to S_{\hz}(r)$.  If we define $F_{\hz}=Exp_{\hz}^{-1}(V \cap \partial \hm)$, there exists $r^{\prime} < r$ such that $F_{\hz} \subset V_{\hz} \setminus {\overline B}_{\hz}(r^{\prime})$.  Since $\pi : V_{\hz} \setminus {\overline B}_{\hz}(r^{\prime}) \to S_{\hz}(r)$ is a Lipschitz map. Thus $\pi(F_{\hz})$ has Hausdorff dimension $< m-1$ in $S_{\hz}(r)$ (which is $(m-1)$-dimensional). Its complementary $F_{\hz}^c$ is thus a dense open subset of $S_{\hz}(r)$.  Also $Exp_{\hz}(\pi^{-1}(F_{\hz}^c))$ is a dense open  subset of $V \cap \hm$, which is included in $W_{\hz}$. This proves that $W_{\hz}$ itself is dense. It's not difficult to check that $W_{\hz}$ is also open.  

\end{preuve}

Thanks to this lemma, $W_{\hx_{k_0}}$ is dense in $V_{\hx} \cap \hm$.  This later set is itself dense in $V_{\hx} \cap (\hl \cup \hm)$.  By continuity of $\overline \sigma$, we get that $\overline \sigma$ and $\sigma_{\hx}^{\prime}$ coincide on $V_{\hx} \cap (\hl \cup \hm)$, as desired.

Assume now that $\overline \Omega$ is a hypersurface with boundary, and $\partial \Omega$ is locally  the graph of a Lipschitz function. Let $X \in \lieg$ and $X_{\Lg}$ (resp. $X_{N}$) be the smooth vector field  on $\hlg$ (resp. on $\hn$) defined by $\omega^{\Lg}(X_{\Lg})=X$ (resp. $\omega^{N}(X_{N})=X$).  The local flows generated by $X_{\Lg}$ and $X_N$ are denoted by $\phi_{X_{\Lg}}^t$ and $\phi_{X_{N}}^t$ respectively. Let $W$ be a sm�all domain   of  a topological hypersurface in $\hlg$ (resp. in $\hn$) which is transverse to $X_{\Lg}$ (resp. to $X_N$) on the interval $]-\epsilon^{\prime},\epsilon[$. We mean here that for $t \not = t^{\prime}$ in $]-\epsilon^{\prime},\epsilon[$, we have $\phi_{X_{\Lg}}^t(W) \cap \phi_{X_{\Lg}}^{t^{\prime}}(W)= \emptyset$ (resp. $\phi_{X_{N}}^t(W) \cap \phi_{X_{N}}^{t^{\prime}}(W)= \emptyset$). We then define:

$$ B_{W,X_{\Lg}}^{]-\epsilon^{\prime},\epsilon[}=\bigcup_{t \in ]-\epsilon^{\prime},\epsilon[}\phi_{X_{\Lg}}^t(W)$$
and
$$  B_{W,X_{N}}^{]-\epsilon^{\prime},\epsilon[}=\bigcup_{t \in ]-\epsilon^{\prime},\epsilon[}\phi_{X_{N}}^t(W)$$

Those are open neighbourhoods in $\hlg$ and $\hn$ respectively. We saw in the proof of Lemma \ref{identification} that $\hl$ is the projection by the smooth covering map $\pi_{\Lg} : \hL \to \hlg$ of the inverse image by $\pi_L : \hL \to L$ of $\tilde \Lambda$, an open subset of $\partial \Omega$. We thus see that  $\hl$ is  locally the graph of a Lipschitz function. As a consequence, given $\hx \in \hl$, there exists $W$ a neighbourhood of $\hx$ in $\hl$, and a vector $X \in \lieg$ such that $X_{\Lg}$ is transverse to $W$, i.e there exists $\epsilon >0$ s.t $\phi_{X_{\Lg}}^t(W) \cap \phi_{X_{\Lg}}^{t^{\prime}}(W)= \emptyset$ for $t \not = t^{\prime}$ in $]-\epsilon,\epsilon[$.  Considering, if necessary,  $-X$ instead of $X$, we will assume that $\phi_{X_{\Lg}}^t.\hx \in \hm$ for $t \in ]- \epsilon, 0[$. Let   $B$ be a neighbourhood of $\hy=\partial \sigma(\hx)$ in $\hn$, which is moreover a flow-box for $\phi_{X_N}^t$, and with the property that for each $\hz \in B$, $\phi_{X_N}^t.\hz$ is defined on $[-\epsilon_0,+\epsilon_0]$.    Choosing $\epsilon < \epsilon_0$ and $W$ small enough, we may assume that $\phi_{X_{\Lg}}^{-\frac{\epsilon}{2}} (W) \subset \hm$ and that $\sigma \circ \phi_{X_{\Lg}}(W) \subset B$.  

Let us call $W^{\prime}=\phi_{X_{\Lg}}^{-\frac{\epsilon}{2}} (W)$. We choose $W^{\prime \prime}$ a small piece of  smooth hypersurface, which is contained in $B_{W^{\prime},X_{\Lg}}^{]-\epsilon/2,0[}$ and transverse to $X_{L_{\Gamma}}$, and $\epsilon_0>\epsilon^{\prime}>0$ such that $V_{\hx}=B_{W^{\prime \prime},X_{\Lg}}^{]-\epsilon^{\prime},\epsilon[}$ is an open neighbourhood of $\hx$. Every point of $V_{\hx}$ writes in an unique way $\hz=\phi_{X_{\Lg}}^{t_z}.{\hat w_z}^{\prime \prime}$, with ${\hat w_z}^{\prime \prime} \in W^{\prime \prime}$, and $t_z\in ]-\epsilon^{\prime},\epsilon[$. 

We then define a map $\sigma_{\hx}^{\prime}$ on $V_{\hx}$ by the formula:
$$ \sigma_{\hx}^{\prime}(\hz)=\phi_{X_N}^{t_z} \circ \sigma ({\hat w_z}^{\prime \prime})$$
Since $\hz \mapsto t_z$ and $\hz \mapsto w_z^{\prime \prime}$ are smooth, $\sigma_{\hx}^{\prime}$ is a smooth map defined on $V_{\hx}$.  Now,  $V_{\hy}=B_{\sigma(W^{\prime \prime}),X_{N}}^{]-\epsilon^{\prime},\epsilon[}$  is an open neighbourhood of  $\hy$ in $\hn$. Since it is included in $B$, it is a flow-box for $\phi_{X_N}^t$, what implies that $\sigma_{\hx}^{\prime}$ is  a smooth diffeomorphism between $V_{\hx}$ and $V_{\hy}$.

The fact that $\sigma$ is a geometrical embedding yields:

$$\sigma \circ \phi_{X_{\Lg}}^t (\hx) = \phi_{X_N}^t \circ \sigma(\hx)$$
 for every $t \in \re$ and $\hx \in \hm$ such that $\phi_{X_{\Lg}}^t (\hx) \in \hm$ and $\phi_{X_{N}}^t \circ \sigma (\hx)$ is  defined. In particular $\sigma_{\hx}^{\prime}$ coincides with $\sigma$ on $V_{\hx} \cap \hm$, and with $\overline \sigma$ on $V_{\hx} \cap (\hl \cup \hm)$ by continuity of $\overline \sigma$.
\end{preuve}



\section{Rigidity of conformal embeddings: the Riemannian case}
 \label{riemannian-case}
 We are going to apply the general tools of the previous sections to the particular case of Riemannian conformal embeddings. Before this, let us recall why conformal Riemannian structures can be described as Cartan geometries infinitesimally modelled on the sphere ${\bf S}^n=PO(1,n+1)/P$.  

\subsection{Conformal structures from the point of view of Cartan connections}
Let us denote by $O(1,n+1)$ the subgroup of $GL(n+2,\re)$ preserving the quadratic form $-x_1^2+x_2^2+...+x_{n+2}^2$, and by $\lieo(1,n+1)$ its Lie algebra.  Recall that $PO(1,n+1):=O(1,n+1)/\{  \pm Id\}$ coincides with the group of conformal transformations of ${\bf S}^n$ endowed with the round metric. Thus, we see the conformal sphere as the homogeneous space ${\bf S}^n=PO(1,n+1)/P$, where $P$ is the stabilizer of an isotropic line for $-x_1^2+x_2^2+...+x_{n+2}^2$.  Using a stereographic projection, $P$ can be seen as the conformal group of the Euclidean space, hence is isomorphic to the semi-direct product $(\re \times O(n)) \ltimes \re^n$. The Lie algebra    $\lieo(1,n+1)$ is a sum $\lien^+ \oplus {\bf R} \oplus \lieo(n) \oplus \lien^-$, where $\liep= {\bf R} \oplus \lieo(n) \oplus \lien^-$ is the Lie algebra of $P$. The algebras $\lien^+$ and $\lien^-$ are both abelian of dimension $n$. 

Assume now that $(M,\hm,\om)$ is a Cartan geometry modelled on ${\bf S}^n=PO(1,n+1)/P$.  For $x \in M$, each $\hx \in \hm$ above $x$ determines an isomorphism $i_{\hx} : T_xM \to \lieo(1,n+1)/\liep$ in the following way: if $u \in T_xM$ and $\hat u$ is a lift of $u$ in $T_{\hx}\hm$, $i_{\hx}(u)$ is the projection on  $\lieo(1,n+1)/\liep$ of $\omega_{\hx}^M(\hat u)$ (this does not depend on the lift $\hat u$). Since $R_p^* \omega = \Ad p^{-1} \omega$, for any $p \in P$, we get:

\begin{equation}
\label{equi}
i_{\hx.p}=\Ad p^{-1} i_{\hx}
\end{equation} 
On $\lieo(1,n+1)$, there is a unique $\Ad P$-invariant conformal class ${\cal C}$ of Riemannian scalar products.  Pulling back this class by $i_{\hx}$, we endow each $T_xM$ with a conformal class of Riemannian scalar products. The construction does not depend on the choice of $\hx$ above $x$ by the identity (\ref{equi}). Hence, any data $(M,\hm,\om)$ as above defines a conformal class of Riemannian metrics on $M$.  

What is nice is that this construction can be reversed. Indeed, it is known since Elie Cartan that given a conformal class of Riemannian metrics on a manifold $M$ of dimension $\geq 3$, there is on the  bundle $\hm$ of $2$-jets of orthogonal frames a {\it unique} normal Cartan connection $\om$ with values in $\lieo(1,n+1)$.  The normality condition is put on the curvature of the connection to ensure uniqueness (see \cite{kobayashi}, \cite{sharpe}). Of course, the normal Cartan connection gives back the initial conformal class by the construction decribed above.

Now, if $(M,g)$ and $(N,h)$ are two Riemannian manifolds of the same dimension $n \geq 3$, and $s : (M,g) \to (N,h)$ a conformal embedding, and if $(M,\hm,\om)$ and $(N,\hn,\on)$ are the normal Cartan geometries modelled on $PO(1,n+1)/P$ canonically associated, then $s$ lifts to $\sigma : (M,\hm,\om) \to (N,\hn,\on)$ a geometrical embedding.

\subsection{Conformal cones}
Let $(L,{\hat L},\ol)$ be a Cartan geometry modelled on ${\bf S}^n=PO(1,n+1)/P$.  For $\hx \in \hL$ and $u \in T_{\hx}\hL$, we can consider the vector field $U_L$ on $\hL$ such that $\ol(U_L)=\omega_{\hx}^L(u)$, and the associated local flow $\phi_{U_L}^t$. We can now introduce the exponential map $\exp_{\hx}$ at $\hx$, defined from a neighbourhood $W_{\hx}$ of $0_{\hx}$ in $T_{\hx}\hL$ to $\hL$. By definition $\exp_{\hx}(u)=\phi_{U_L}^1.\hx$.  Notice that $\exp_{\hx}$ is a diffeomorphism from a  sufficiently small neighbourhood of $0_{\hx}$ on its image.  

We are now going to define what is a {\it conformal cone} on a Riemannian manifold $(L,g)$.  By what has been said before, $(L,g)$ defines a unique normal Cartan geometry $(L,\hL,\ol)$ modelled on ${\bf S}^n=PO(1,n+1)/P$.  The subalgebra $\lien^+ \subset \lieo(1,n+1)$ is not $\Ad P$-invariant.  Nevertheless, writting $P=(\re \times O(n)) \ltimes \re^n$, we see that $\lien^+$ is invariant by the adjoint action of $\re \times O(n)$. We put on $\lien^+$ a scalar product $< >$ which is $\Ad O(n)$-invariant, and call $||.||$ the asociated norm.  Let $S_{\lien^+}$ be the unit sphere of $\lien^+$ for the norm $||.||$, and ${\cal B}$ a closed ball of nonzero radius of $S_{\lien^+}$, for the metric induced by $<>$ on $S_{\lien^+}$. Given $\lambda >0$, we can then define:
$$ {\cal C}({\cal B}, \lambda) = \{ u \in \lien^+ \ | \ u=tv \  t \in [0,\lambda], \   v \in {\cal B}   \}$$  
If $\pi_{\X}$ denotes the projection $PO(1,n+1) \to PO(1,n+1)/P$, then we define a subset of ${\bf S}^n$ putting $C({\cal B},\lambda)= \pi_{\X} \circ \exp_G({\cal C}({\cal B},\lambda))$.
Calling $o$ the fixed point of $P$ on ${\bf S}^n$, such a subset will be referred to as a cone with vertex $o$.

Now, for any $x \in L$, and $\hx \in \hL$ above $x$, we can define ${\cal C}_{\hx}({\cal B},\lambda)=(\omega_{\hx}^L)^{-1}({\cal C}({\cal B},\lambda))$, and ${\hat C}_{\hx}({\cal B},\lambda)=exp_{\hx}({\cal C}_{\hx}({\cal B},\lambda))$. Let $\pi : \hL \to L$ be the bundle projection. A set of the form $C_{\hx}({\cal B}, \lambda)=\pi({\hat C}_{\hx}({\cal B}, \lambda))$  will be called a {\it conformal cone with vertex $x$}. We will also use the notation ${\dot C}_{\hx}({\cal B}, \lambda)$, which is simply ${C}_{\hx}({\cal B}, \lambda)$ with its vertex $x$ removed.  Since a conformal transformation of $(L,g)$ lifts to an automorphism of $(\hL,\ol)$, it is not difficult to check that the set of conformal cones is preserved by conformal transformations. 

Also useful for what follows will be the notion of {\it development} of a conformal cone. Let $x \in L$, $\hx \in \hL$ above $x$, $\hx^{\prime}=\hx.p$ with $p \in P$, and $C_{\hx^{\prime}}({\cal B},\lambda)$ a conformal cone with vertex $x$. We set:
$$ {\cal D}_{x}^{\hx}(C_{\hx^{\prime}}({\cal B},\lambda)):=p.C({\cal B},\lambda)$$

We won't detail here the notion of development of curves for a Cartan geometry, which is of course very close to the previous definition (see \cite{sharpe} chapter $5$ for more details).  

\subsection{The boundary map for Riemannian conformal embeddings}
\label{boundary-riemannian}
Until the begining of section \ref{kleinian-manifold}, we assume that $(L,g)$ is a smooth Riemannian manifold of dimension $n \geq 3$, and $M \subsetneq L$ is an open subset. We assume that we are in one of the following cases:
\begin{enumerate}
\item{The Hausdorff dimension of $\partial M$ is $<n-1$.}
\item{The closure $\overline M$ is a topological manifold with locally Lipschitz boundary of $L$ (see hypothesis ${\bf H_1}$ of the introduction).}
\end{enumerate}
We assume that $s : (M,g) \to (N,h)$ is a {\it strict} conformal embedding, where $(N,h)$ is a Riemannian manifold of dimension $n$. This conformal embedding lifts to a geometrical embedding  $\sigma :(M,\hm,\om) \to (N, \hn, \on)$
between the normal Cartan geometries defined by the two conformal structures. As in section  \ref{natural-map}, we introduce the set:
$$ \hl = \{  \hx \in \partial \hm \ | \ \exists (\hx_k) \ a \ sequence \ of \ \hm, \ \hx_k \to \hx \ s.t \ \sigma(\hx_k) \ converges \ in \ \hn  \}$$
and $\Lambda=\hl/P$.
It turns out that in the case of Riemannian conformal structures, we understand very well $\Lambda$ and the boundary map $\partial s : \Lambda \to \partial_sM$.  The key proposition is:

\begin{proposition}
\label{property-H}
Let $(L,g)$ and $M \subsetneq L$ satisfying one of the conditions $1.$ or $2.$ above. Let $s : (M,g) \to (N,h)$ be a strict conformal embedding.  Then: 
\begin{enumerate}
\item{ $\Lambda$ coincides with the set:
$$ \{ x \in \partial M \ | \ \exists (x_k) \ a\ sequence \ of \ M, \ x_k \to x\ and \ s(x_k) \ converges  \ in \ N \}$$}
\item{If  $\overline M$ is compact, then $\Lambda $ is a nonempty open subset of $\partial M$, and the map $\partial s: \Lambda \to \partial_sM$  is surjective and proper.}

\end{enumerate}

\end{proposition}

\begin{preuve}
By the very definition of $\Lambda$,  one always have the inclusion $\Lambda \subset \{ x \in \partial M \ | \ \exists (x_k) \ a\ sequence \ of \ M, \ x_k \to x\ and \ s(x_k) \ converges  \ in \ N \}$. If this latter set is empty, so is $\Lambda$. 

Let $x \in \partial M$. Assume that there is a sequence $(x_k)$ of $M$ tending to $x$, such that $s (x_k)=y_k$ tends to $y \in \partial_sM$. 
We choose a sequence $(\hx_k)$ in $\hm$ (projecting onto $(x_k)$), converging to $\hx \in \partial \hm$, and a sequence $(p_k)$ in $P$ such that $\hy_k=\sigma(\hx_k).p_k^{-1}$ tends to $\hy \in \partial_{\sigma}\hm$.  Our aim is to prove that $\hx \in \hl$, which is easily seen to be equivalent to the sequence $(p_k)$ being bounded in $P$.

Assume by contradiction that $(p_k)$ is unbounded. We will take advantage of the assumptions made on $\partial M$  thanks to  the:

\begin{lemme}
\label{exemples-propriete-H}
If $\hm \subset \hat L$ is an open subset  such that either the Hausdorff dimension of $\partial \hm$ is $< \text{dim} \hat L-1$, or $\partial \hm$ is locally Lipschitz, then for every $\hx \in \partial \hm$ there is a cone $\hat C_{\hx}({\cal B},\lambda)$ such that $\hat C_{\hx}({\cal B},\lambda) \setminus \{ \hx \} \subset \hm$. 
\end{lemme}
Let us choose such a cone $\hat C_{\hx}({\cal B},\lambda)$. Taking ${\cal B}$ and $\lambda$ smaller if necessary, we can assume that there is a $k_0 \in \NN$ such that $\hat C_{\hx_k}({\cal B},\lambda) $ is in the interior of $\hat C_{\hx}({\cal B},\lambda)$ as soon as $k \geq k_0$, hence  $\hat C_{\hx_k}({\cal B},\lambda) \subset \hm$.  We have $s(C_{\hx_k}({\cal B}^{\prime}, \lambda))=C_{\sigma(\hx_k)}({\cal B}^{\prime}, \lambda)$, so that:
\begin{equation}
\label{equation-cones}
 {\cal D}_{y_k}^{\hy_k}(C_{\sigma(\hx_k)}({\cal B}^{\prime},\lambda))=p_k.C({\cal B}^{\prime},\lambda)
 \end{equation} 
We are now confronted to the following problem: first, determine the behaviour of the sequence of sets $p_k.C({\cal B}^{\prime},\lambda)$ in ${\bf S}^n$, when $(p_k)$ is a sequence of $P$ tending to infinity. When it will be done, consider the natural question: {\it given a sequence of conformal cones on $M$, can we deduce the behaviour of this sequence (for the Hausdorff topology), knowing the behaviour of its developments}?  The two following lemmas answer two our question.

\begin{lemme}
\label{dynamique-cones}
Let $C({\cal B},\lambda)$ be a conformal cone with vertex $o$, $\lambda >0$. Let $(p_k)$ be a sequence of $P$ tending to infinity. Then, considering a subsequence of $(p_k)$ if necessary, we are in one of the following cases:

$(i)$ there exists ${\cal B}^{\prime} \subset {\cal B}$ a closed subball with nonzero radius, such that $p_k.C({\cal B}^{\prime}, \lambda)\to o$.

$(ii)$ there exists  a cone ${\cal C}({{\cal B}_{\infty}},\lambda_{\infty})$,  a sequence $(\epsilon_k)$ of ${\Bbb R}_*^+$ tending to $0$, a sequence $l_k$ tending to $l_{\infty}$ in $P$ such that $l_kp_k$ stays in the factor $\re \times O(n)$ of $P=(\re \times O(n)) \ltimes \re^n$, and $\Ad (l_kp_k).{\cal C}({\cal B},\epsilon_k) \to {\cal C}({{\cal B}_{\infty}},\lambda_{\infty})$. 
\end{lemme}

\begin{lemme}
\label{to-zero}
Let $(L,[g])$ be a Riemannian conformal structure of dimension $\geq 3$. Let $(x_k)$ be a sequence of $L$ converging to $x$, and $(\hx_k)$, $(\hx_k)^{\prime}$ two sequences of $\hl$, projecting on $(x_k)$. We assume that $(\hx_k)$ converges to $\hx \in \hl$.  Then, if ${\cal D}_{x_k}^{\hx_k}(C_{\hx_k^{\prime}}({\cal B},\lambda)) \to o$ for the Hausdorff topology in ${\bf S}^n$, then $C_{\hx_k^{\prime}}({\cal B},\lambda) \to x$ for the Hausdorff topology on $L$.   

\end{lemme}

We postpone the proof of lemma \ref{dynamique-cones} to section \ref{preuve-dynamique}. Lemma \ref{to-zero} is a particular case of \cite{frances1}, Lemma 7.

We now use the conclusions of Lemma \ref{dynamique-cones}:

$\bullet$ either there exists ${\cal B}^{\prime} \subset {\cal B}$, such that $p_k.C({\cal B}^{\prime},\lambda) \to o$.  Then, by the relation (\ref{equation-cones}) and Lemma \ref{to-zero}, we conclude  that $s(C_{\hx_k}({\cal B}^{\prime}, \lambda)) \to y$ for the Hausdorff topology. But this is a contradiction. Indeed, since for any $\hx \in \hl$, $\exp_{\hx}$ is a diffeomorphism from a  a neighbourhood $U_{\hx}$ of $0_{\hx}$ on its image, and since $\exp_{\hx}(U_{\hx} \cap (\omega_{\hx}^{L})^{-1}(\lien^+))$ is transverse to the fibers of $\pi : \hat L \to L$ for $U_{\hx}$ sufficiently small, we deduce that any conformal cone $C_{\hx}({\cal B},\lambda)$ has nonempty interior. From this, it follows that all the sets $C_{\hx_k}({\cal B}^{\prime}, \lambda)$, $k \geq k_0$, contain a common open subset $U \subset M$. We thus get: $\sigma(U) = \{ y \}$, which yields a contradiction.

$\bullet$ if we are in the second case of Lemma \ref{dynamique-cones}, then:
$$ \Ad l_kp_k. {\cal C}({\cal B}, \epsilon_k) \to {\cal C}({\cal B}_{\infty}, \lambda_{\infty})$$
where ${\cal C}({\cal B}_{\infty}, \lambda_{\infty})$ is a cone of $\lien^+$,  $\epsilon_k \to 0$, $l_k \to l_{\infty}$ in $P$ and $l_kp_k$ is in  $\re \times O(n) \subset P$. In particular $\Ad l_kp_k. {\cal C}({\cal B}, \epsilon_k)$ is a sequence of cones of $\lien^+$ containing a common open subset ${\cal U}$ for $k$ large. Now, putting $\hy_k^{\prime}=\hy_k.l_k^{-1}$, we get that $\hy_k^{\prime} \to \hy.l_{\infty}^{-1}$. Also, 
$$ \sigma({\hat C}_{\hx_k}({\cal B}, \epsilon_k))={\hat C}_{\sigma(\hx_k)}({\cal B}, \epsilon_k)$$
what means that 
$$\sigma(\hat C_{\hx_k}({\cal B}, \epsilon_k))=exp_{\hy_k^{\prime}}(\Ad l_kp_k. {\cal C}_{\hy_k^{\prime}}({\cal B}, \epsilon_k))$$
Now $exp_{\hy_k^{\prime}}(\Ad l_kp_k. {\cal C}_{\hy_k^{\prime}}({\cal B}, \epsilon_k))$ contains $exp_{\hy_k^{\prime}}((\omega_{\hy_k^{\prime}}^L)^{-1}({\cal U}))$ for $k$ sufficiently large. If ${\cal U}$ was chosen close to $0_{\lien^+}$, then $exp_{\hy_k^{\prime}}((\omega_{\hy_k^{\prime}}^L)^{-1}({\cal U}))$ is transverse to the fibers of $\pi_N : \hn \to N$, and $\pi_N \circ exp_{\hy_k^{\prime}}((\omega_{\hy_k^{\prime}}^L)^{-1}({\cal U}))$ is an open subset of $s(M)$. Since $\hy_k^{\prime} \to \hy.l_{\infty}^{-1}$, all the open subsets $\pi_N \circ exp_{\hy_k^{\prime}}((\omega_{\hy_k^{\prime}}^L)^{-1}({\cal U}))$ contain a common open set $U \subset s(M)$ for $k$ large.  On the one hand, we thus have that $C_{\hx_k}({\cal B}, \epsilon_k)$ tends to $x$ for the Hausdorff topology, hence leaves every compact subset of $M$, and on the other hand the sets  $s(C_{\hx_k}({\cal B}, \epsilon_k))$ contain a common open subset $U \subset s(M)$ for $k$ large. This contradicts that $s : M \to s(M)$ is a proper map, and the first point of the proposition is proved.


Let us summarize what we got before: if for  $x \in \partial M$, there is $(x_k)$ a sequence of $M$ converging to $x$ such that $s(x_k)$  converges to $y \in N$, then $x \in \Lambda$ and $\partial s(x)=y$.  Now, if $y \in \partial_sM$, there is a sequence $(x_k)$ of $M$ such that $s(x_k) \to y$. If we assume that $\overline M$ is compact in $L$, a subsequence of $(x_k)$ will converge to $x \in \partial M$. Hence $x \in \Lambda$ and $\partial s(x)=y$ by what we just said.  This proves that $\partial s : \Lambda \to \partial_sM$ is surjective in this case.  To get the second point of the proposition, we must show that $\partial s$ is moreover proper.   To this end, consider $(x_k)$  a sequence of $\Lambda$ leaving every compact subset of $\Lambda$. Since $\partial M$ is supposed to be compact, we will assume that $x_k \to x$ for some $x \in \partial M$.  Now, if $\partial s(x_k)$ does not leave every compact subset of $N$, we can  suppose $\partial s(x_k) \to y$ for some $y \in N$. Using the continuity of $\partial s$, one can exhibit a sequence $(x_k^{\prime})$ of $M$, also converging to $x$, such that $s(x_k^{\prime})$ tends to $y$. By the first point of the theorem, we get $x \in \Lambda$: a contradiction.

\end{preuve}

\subsubsection{Proof of lemma \ref{dynamique-cones}}
\label{preuve-dynamique}

To complete the proof of proposition \ref{property-H}, we still have to prove lemma \ref{dynamique-cones}.
Using a stereographic projection, we identify conformally ${\bf S}^n \setminus \{ o \}$ with $\re^n$. Any sequence of $P$ becomes a sequence of conformal transformations of the Euclidean space $p_k : x \mapsto \lambda_k A_kx+\mu_k u_k$, where $\lambda_k \in \Bbb R_+^*$,  $\mu_k \in \Bbb R_+$, $A_k \in O(n)$, and $u_k \in S(1)$. Now, looking at a subsequence if necessary, we assume that $\lambda_k$, $\mu_k$, $\frac{\lambda_k}{\mu_k}$ all have limits in ${\Bbb R}_+^* \cup \{  +\infty\}$,  $u_k \to u_{\infty}$, and $A_k \to A_{\infty}$ in $O(n)$. The conclusions of the lemma won't be affected if we replace $p_k$ by $(A_k)^{-1}.p_k$, so that we will assume now that $p_k= \lambda_k Id + \mu_k u_k$.

To understand the dynamics of a sequence $(p_k)$ on the set of conformal cones with vertex $o$, it is better if we describe more precisely  such a cone when seen through the stereographic projection. The map $s^+ : \lien^+ \setminus \{0 \} \to \re^n \setminus \{ 0 \}$ is a conformal chart, which maps lines of $\lien^+$ through zero to lines of $\re^n$ through $0$. For a suitable choice of $< >$, $s^+$ sends $S_{\lien^+}$, the sphere of center $0_{\lien^+}$ and radius $\lambda >0$ (for $|| . ||$) to $S(\frac{1}{\lambda})$, the Euclidean sphere of radius $\frac{1}{\lambda}$ and center $0$. Thus, in the chart $s^+$, a cone ${\cal C}({\cal B}, \lambda)$ (with the origin removed) just read as the set:
$$ \{ u \in \re^n \ | \ \exists t \in [\frac{1}{\lambda};+ \infty[, \ \frac{u}{t} \in s^+({\cal B})  \}$$

 The following lemma, which proof is easy will also be useful:

\begin{lemme}
\label{dynamique-segments}
Let $[x_k,u_k)$ be a sequence of half-lines in $\re^n$. Assume that there are $u_{\infty}$ and $v_{\infty}$ in $S(1)$ such that $u_k \to u_{\infty}$ and $\frac{x_k}{||x_k||} \to v_{\infty}$.  If $x_k \to \infty$ and $v_{\infty} \not = u_{\infty}$, then $[x_k,u_k) \to [o]$.

\end{lemme}

Assume first that $\mu_k$ tends to $a \in {\Bbb R}_+$. Then, we call $l_k$ the translation of vector $-\mu_ku_k$. Clearly, $l_k \to l_{\infty}$ in $P$, whith $l_{\infty}$ the translation of vector $-au_{\infty}$, and $l_kp_k$ is just the homothetic transformation $x \mapsto \lambda_k x$, hence is in $\re \times O(n) \subset P$. Since $(p_k) \to \infty$, we must have $\lambda_k \to + \infty$ or $\lambda_k \to 0$. In the first case, $p_k.C({\cal B},\lambda) \to [o]$, so that we are in case $(i)$ of the lemma. If $\lambda_k \to 0$, then $l_kp_k.C({\cal B},\lambda_k) \to C({\cal B},1)$ and also $\Ad l_kp_k.{\cal C}({\cal B},\lambda_k) \to {\cal C}({\cal B},1)$, so that we are in case $(ii)$ of the lemma.  

If $\mu_k \to + \infty$, and $\frac{\lambda_k}{\mu_k} \to b$, $b \in {\Bbb R}_+$. Then, let ${\cal B}^{\prime} \subset {\cal B}$ be a closed subball with nonzero radius, such that $- u_{\infty} \not \in {\cal B}^{\prime}$.  Then, there exists $\alpha_0>\frac{1}{\lambda}$, such that if $x \in {\cal B}^{\prime}$, and $\alpha > \alpha_0$, then $\frac{b \alpha x + u_{\infty}}{||�b \alpha x + u_{\infty}||} \not = -u_{\infty}$.  Now, if $[x_k,u_k)$ is a sequence in $C({{\cal B}^{\prime}},\frac{1}{2 \alpha_0})$, Lemma \ref{dynamique-segments} ensures that $p_k.[x_k,u_k) \to [o]$. We are thus in the situation $(i)$ of the lemma.

It remains to investigate the case where $\mu_k \to + \infty$ and $\frac{\lambda_k}{\mu_k} \to + \infty$. Let ${\cal B}^{\prime} \subset {\cal B}$ be a closed subball with nonzero radius, such that $-u_{\infty} \not \in {\cal B}^{\prime}$. Let $[x_k,u_k)$ be a sequence of $C({{\cal B}^{\prime}},\lambda)$.  Writting $p_k.[x_k,u_k)=[x_k^{\prime},u_k)$, with $u_k^{\prime}=\lambda_k(x_k +\frac{\mu_k}{\lambda_k}u_k)$, we get by Lemma \ref{dynamique-segments} that $p_k.[x_k,u_k) \to [o]$. Since this is true for any sequence, we get $p_k.C({{\cal B}^{\prime}},\lambda) \to [o]$. We are once again in the case $(i)$ of the lemma.

\subsection{Embeddings of open subsets with ``small" boundary:  proof of Theorem \ref{codimension2}}
We consider here $(L,g)$ a {\it compact} Riemannian manifold of dimension $n \geq 3$,  $M \subsetneq L$ an open subset such that the boundary $\partial M \subset L$ has Hausdorff dimension $<n-1$, and $s : (M,g) \to (N,h)$  a strict, smooth conformal embedding. By proposition \ref{property-H}, $\Lambda$ is a nonempty open subset of $\partial M$. The hypothesis made on the Hausdorff dimension of $\partial M$ ensures that $\Lambda \cup M$ is actually a dense open subset  $M^{\prime} \subset L$ containing $M$.  Its inverse image $\hm^{\prime} \subset \hL$ is also a dense open subset.  Calling $\omega^{M^{\prime}}$ the restriction of $\ol$ to $\hm^{\prime}$, we get a Cartan geometry $(M^{\prime},\hm^{\prime},\omega^{\prime})$ modelled on ${\bf S}^n$, which is just the canonical Cartan geometry associated to the conformal structure induced by $g$ on $M^{\prime}$.  By proposition \ref{boundary}, we know that $\sigma$ extends to a continuous map $\overline \sigma : \hm^{\prime} \to \hn$ (recall that $\overline \sigma$ coincides with $\partial \sigma$ on $\partial \hm$ and with $\sigma$ on $\hm$).  For every $\hx \in \hl$, the open subset $V_{\hx}$ given by Proposition \ref{extension} can be chosen arbitrarely small, so that we can assume $V_{\hx} \subset \hm^{\prime}$. As a consequence, $\overline \sigma$ coincides with $\sigma_{\hx}^{\prime}$ (given by Proposition \ref{extension}) on $V_{\hx}$, so that $\overline \sigma$ is actually smooth.  The identity $\overline \sigma^* \omega^N=\omega^{M^{\prime}}$ holds on $\hm$, because $s$ is a conformal map.  But $\hm$ is dense in $\hm^{\prime}$, so that  $\overline \sigma^* \omega^N=\omega^{M^{\prime}}$ on $\hm^{\prime}$.  Let us also remark that since $\hm$ is dense in $\hm^{\prime}$ and $\sigma$ is one-to-one, then $\overline \sigma$ is also one-to-one. The map $\overline \sigma$ induces a conformal embedding $\overline s : (M^{\prime},g) \to (N,h)$.  Now, by the second point of Proposition \ref{property-H}, we know that since $\overline M=L$ is compact, $\partial s(\Lambda)=\partial_sM$. As a consequence, the Hausdorff dimension of $\partial_sM$ is $<n-1$, what means that $s(M)$ is dense in $N$ and $\overline s(M^{\prime})=\partial_sM \cup s(M)=N$.  Thus, $\overline s$ is onto and the theorem is proved.

\subsection{Embeddings of open subset bounded by hypersurfaces: Theorems \ref{immersion} and \ref{geometrique}}
\subsubsection{Proof of Theorem \ref{immersion}}
\label{preuve-immersion}
We assume here that $M \subsetneq L$ is an open subset of the Riemannian manifold $(L,g)$ of dimension $n \geq 3$.  The closure $\overline M$ is assumed to be a  compact topological submanifold with boundary of $L$. This boundary  $\partial M$ is a closed ${\cal C}^{k,\alpha}$-hypersurface of $L$, $k+\alpha \geq 1$. We consider a strict conformal embedding $s: (M,g) \to (N,h)$.  

By Lemma \ref{identification} and point $2.$ of Proposition \ref{property-H}, the set $\Lambda$ is nonempty and open in $\partial M$, and is {\it mapped surjectively} on $\partial_sM$.  Now, by Proposition \ref{extension}, each $\hx \in \hl$ has a neighbourhood $V_{\hx}$ such that $\partial \sigma_{| V_{\hx} \cap \hl}$ is the restriction of a smooth diffeomorphism $\sigma_{\hx}^{\prime} : V_{\hx} \to V_{\hy}$ (with $\hy=\partial \sigma(\hx)$). It follows that $\partial \sigma : \hl \to \hn$ is a $P$-equivariant ${\cal C}^{k,\alpha}$-immersion. Hence, $\partial s : \Lambda \to N$ is also a ${\cal C}^{k,\alpha}$ immersion. 

We now finish the proof of the first point of theorem \ref{immersion}, assumin that $k \geq 1$.
\begin{lemme}
\label{lemma-immersion}
If $k \geq 1$, the map $\partial s : \Lambda \to N$ is a conformal immersion, whose image is $\partial_sM$.
\end{lemme}

\begin{proof}  
 We already know that $\partial s$ is a ${\cal C}^k$-immersion, which is onto $\partial_sM$, and we have to prove that it is conformal.  To see this, let us recall that to each $\hx \in \hat L$ over $x$, we can associate  a natural isomorphism $i_{\hx} : T_xL \to \lieo(1,n+1)/\liep$.  It is built in the following way: for $u \in T_xL$, choose $\hat u \in T_{\hx}\hat L$ such that $\pi_*(\hat u)=u$, and define $i_{\hx}(u):=\omega_{\hx}^L(\hat u)$.  It is clear that this definition does not depend on the choice of $\hat u$. By the very definition of the normal Cartan connection $\omega^L$ associated to the conformal structure $(L,[g])$, if $<>$ is a scalar product on $\lieo(1,n+1)/\liep$ which is $\Ad O(n)$-invariant, then $(i_{\hx})^*<>$ is in the conformal class $[g]_x$ for every $\hx$ over $x$ (see for example \cite{sharpe}).  If $\hx \in {\hat \Lambda}$, $i_{\hx}$ yields by restriction an embedding from $T_x\Lambda$ into $\lieo(1,n+1)/\liep$.  In a similar way, we have an embedding $i_{\partial \sigma(\hx)} : T_{s(x)}(\partial_sM) \to \lieo(1,n+1)/\liep$.  Since $(\sigma_{\hx}^{\prime})^*\on=\omega^L$ is true on $V_{\hx} \cap \hm$, it remains true on $V_{\hx} \cap (\hm \cup \hat \Lambda)$.  In particular, we get that $i_{\partial \sigma(\hx)} \circ D_x(\partial s)=i_{\hx}$.  This yields:
$$ (D_x(\partial s))^* \circ (i_{\partial \sigma(\hx)})^*<>=(i_{\hx})^*<>$$
what means exactly that $D_x(\partial s)$ sends the restriction of $[g]_x$ to $T_x \Lambda$ to the restriction of $[h]_{\partial s(x)}$ to $T_{\partial s(x)}(\partial_sM)$, i.e $\partial s$ is a conformal map.
\end{proof}

We now prove:
\begin{lemme}
\label{fibers}
If $k \geq 1$, the fibers of $\partial s : \Lambda \to \partial_sM$ have at most two elements.
\end{lemme}
\begin{proof}  Assume that  $x$,$x^{\prime}$ and $x^{\prime \prime}$ are three distinct points of $\Lambda$, which are mapped by $\partial s$ to a same point $y$. We get three points $\hx$, $\hx^{\prime}$ and $\hx^{\prime \prime}$, such that no two of them are in the same fiber, such that $\partial \sigma(\hx)=\hy$, $\partial \sigma(\hx^{\prime})=\hy.p^{\prime}$  and $\partial \sigma(\hx^{\prime \prime})=\hy.p^{\prime \prime}$, for $p^{\prime}$ and $p^{\prime \prime}$ in $P$.   Let us consider $V_{\hx}$,   $V_{\hx}^{\prime}$, $V_{\hx}^{\prime \prime}$ the neighbourhoods given by Proposition \ref{extension}, and $\sigma_{\hx}^{\prime}$, $\sigma_{\hx^{\prime}}^{\prime}$, $\sigma_{\hx^{\prime \prime}}^{\prime}$ the associated diffeomorphisms introduced in the same proposition. Since $\phm$ is a $C^1$ hypersurface, we may assume, shrinking $V_{\hx}$,   $V_{\hx}^{\prime}$, $V_{\hx}^{\prime \prime}$ if necessary, that  $V_{\hx} \cap \hm$, $V_{\hx^{\prime}} \cap \hm$, $V_{\hx^{\prime \prime}} \cap \hm$ are connected and project onto three disjoints open sets $U_x$, $U_{x^{\prime}}$, $U_{x^{\prime \prime}}$ of $M$.  The open sets $\sigma_{\hx}^{\prime}(V_{\hx})$,$\sigma_{\hx^{\prime}}^{\prime}(V_{\hx}^{\prime})$, $\sigma_{\hx^{\prime \prime}}^{\prime}(V_{\hx}^{\prime \prime})$ project on three neighbourhoods of $y$, that we call  $U_y$, $U_y^{\prime}$ and $U_y^{\prime \prime}$.  The projections on $M$ of $\partial \sigma(V_{\hx} \cap \phm)$, $\partial \sigma(V_{\hx^{\prime}} \cap \phm)$ and $\partial \sigma(V_{\hx^{\prime \prime}} \cap \phm)$ are three embedded ${\cal C}^k$-hypersurfaces $\Sigma_y$, $\Sigma_y^{\prime}$ and $\Sigma_y^{\prime \prime}$. By $P$-equivariance of $\overline \sigma$ on $\hl \cup \hm$, the open sets $s(U_x)$, $s(U_{x^{\prime}})$ and $s(U_{x^{\prime \prime}})$ are connected components of $V_y \setminus \Sigma_y$, $V_y^{\prime} \setminus \Sigma_y^{\prime}$ and $V_y^{\prime \prime} \setminus \Sigma_y^{\prime \prime}$ respectively. Hence, $s(U_x)$, $s(U_{x^{\prime}})$ and $s(U_{x^{\prime \prime}})$ can not be pairwise disjoint, contradicting the injectivity of $s$. This shows that the fibers of $\partial s$ have at most two elements.
\end{proof}
The two previous lemmas give the first point of  theorem \ref{immersion}.  To prove the second point, we make the assumption that $s(M)$ has compact closure in $N$.  By the second point of proposition \ref{property-H} which says that $\partial s$ is a proper map, we get that $\Lambda = \partial M$.
We  introduce the following notation: two points $x$ and $x^{\prime}$ on $\Lambda$ are {\it twin points} if $\partial s(x)=\partial s(x^{\prime})$ and $x \not = x^{\prime}$. If $x$ and $x^{\prime}$ are twin points, we say that $x$ branches with $x^{\prime}$ when there exists $U$ a neighbourhood of $x$ in $\Lambda$, such that  for every neighbourhood $U^{\prime}$ of $x^{\prime}$: $\partial s (U^{\prime})  \not \subset \partial s (U)$. Now, define:
$$\Lambda^o = \{  x \in \Lambda \ | \ x \ does \ not \ branch \ with \ any \ x^{\prime} \in \Lambda \} $$

\begin{lemme}
The set $\Lambda^o$ is a dense open subset of $\partial M$.
\end{lemme}

\begin{proof}
Let us call $F$ the complementary of $\Lambda^o$ in $\Lambda$.  If $F$ is empty, there is nothing to prove, so that we make the assumption $F \not = \emptyset$. Assume that $(x_k)$ is a sequence of $F$ converging to $x_{\infty} \in \partial M=\Lambda$. Each point $x_k$ has by definition a twin point $y_k$ in $\partial M$, and considering a subsequence if necessary, we assume that $y_k$  tends to $y_{\infty} \in \Lambda$.  Since by proposition \ref{extension}, $\partial s$ is locally one-to-one, and $\partial s(x_k)=\partial s(y_k)$ for all $k$, we can not have $y_{\infty}=x_{\infty}$.  Thus $y_{\infty}$ is a twin point of $x_{\infty}$.  If $x_{\infty}$ does not branch with $y_{\infty}$, we can find $U$ (resp. $V$) a neighbourhood of $x_{\infty}$ (resp. of $y_{\infty}$) in $\Lambda$ on which $\partial s$ is one-to-one and such that $\partial s (U)= \partial s (V)$. Let $k$ be big enough so that $x_k \in U$ and $y_k \in V$.  Then $\partial s (U)= \partial s (V)$ means that  $x_k$ branches with $y_k$, a contradiction since both of those points are in $F$.  We get that $x_{\infty}$ branches with $y_{\infty}$, what shows that $F$ is closed.

Let us  cover $\Lambda=\partial M$ by a finite family of open subsets $(W_i)_{i=1, ..., s}$ such that for each $i=1,...,s$, the restriction of  $\partial s$ to $W_i$ is one-to-one.  For each $i =1, ..., s $, we introduce the set:
$$ F_i = \{ x \in \Lambda \ | \ x \ branches \ with \ some \ x^{\prime} \in W_i \}$$
We would like to prove that $F_i$ has empty interior.  For a contradiction, assume that some open subset $U \subset \partial M$ is included in $F_i$.  We can assume that $\partial s$ is one-to-one in restriction to $U$.  By definition of $F_i$, every point of $U$ is twin with some point of $W_i$, so that  $\partial s(U) \subset \partial s (W_i)$, and $\partial s (U)$ is open in $\partial s (W_i)$ because $\partial s $ is an immersion.  We would get an open subset $U^{\prime} \subset W_i$ with $\partial s (U^{\prime})= \partial s (U)$, but then, no point of $U$ would branch with any point of $U^{\prime}$: contradiction.

The set $F=\bigcup_{i=1}^sF_i$ must then have empty interior.  We call $\partial_s^oM:=\partial s (\Lambda^o)$.  Let us prove that  $\partial_s^oM$ is a submanifold of $N$.  For that, choose $z \in \partial_s^oM$.  By definition, there is $x \in \Lambda^o$ such that $\partial s (x)=z$.  Let us choose $\Sigma_x$ an open neighbourhood of $x$ in $\partial M$, contained in $\Lambda^o$, and such that $\partial s$ is one-to-one in restriction to $\Sigma_x$. Then,  $\Sigma_z:=\partial s (\Sigma_x)$ is a ${\cal C}^{k,\alpha}$-submanifold of $N$.  We claim that shrinking $\Sigma_x$ (hence $\Sigma_z$) sufficiently, there is an open set $V$ in $N$ such that $V \cap \partial_sM=\Sigma_z$.  If it is not the case, we could find a sequence $(y_k)$ in $\partial M$ such that $z_k:=\partial s (y_k)$ converges to $z$, but $\forall k \in \NN$, $z_k \not \in \Sigma_z$.  Now, we can assume that $y_k$ tends to $y_{\infty}$ in $\partial M$, and because $z_k \not \in \Sigma_z$, $y_{\infty} \not = x$.  Thus $x$ and $y_{\infty}$ are twin points, but since $x \in \Lambda^o$, $x$ does not branch with $y_{\infty}$.  Thus,  we can find $\Sigma_{y_{\infty}}$ a neighbourhood of $y_{\infty}$ in $\partial M$ and $\Sigma_x^{\prime} \subset \Sigma_x$ a neighbourhood of $x$, such that $\partial s (\Sigma_{y_{\infty}})= \partial s(\Sigma_x^{\prime})$. As a consequence, $z_k \in \Sigma_z$ for $k$ big enough: a contradiction.

We just proved that $\partial_s^oM$ is open in $\partial_sM$, and is a ${\cal C}^{k,\alpha}$-submanifold of $N$.  It is dense in $\partial_sM$ since $\Lambda^o$ is dense in $\partial M$ and $\partial s : \partial M \to \partial_s M$ is onto.
\end{proof}

\subsubsection{Proof of Theorem \ref{geometrique}}
\label{preuve-geometrique}
We keep the same hypotheses than for the previous theorem, and we  assume moreover that $\partial_sM$ is a ${\cal C}^0$-hypersurface of $N$.

This hypothesis implies that the set $\Lambda^o$ introduced in the proof of theorem \ref{immersion} is actually equal to $\Lambda$.  Indeed, if $x$ and $x^{\prime}$ are twin points in $\Lambda$ and $x$ branches with $x^{\prime}$, then $\partial s (\Lambda)$ can not be locally the graph of a function.   This implies that $\partial_sM=\partial s(\Lambda)$ is a ${\cal C}^{k,\alpha}$-hypersurface of $N$, and the first point of theorem \ref{geometrique} is proved.  

Assuming moreover that $N$ is  compact, we deduce from point $2.$ of Proposition \ref{property-H} that $\Lambda$ is the whole $\partial M$.  Now, $\partial s : \partial M \to \partial_sM$ is a ${\cal C}^{k,\alpha}$-immersion by theorem \ref{immersion}.  If $\partial M_i$ denotes a connected component of $\partial M$, then $\partial s$ is a $C^{k,\alpha}$-immersion from $\partial M_i$ onto a connected component $\partial_sM_i$ of $\partial_sM$.  Since $\partial M_i$ is compact,  $\partial s$ must be  a covering map from $\partial M_i$ onto $\partial_sM_i$. The fact that this covering is conformal and at most twofold when $k \geq 1$ is a straigthforward consequence of Theorem \ref{immersion}.



\section{Geometric embeddings of Kleinian manifolds}
\label{kleinian-manifold}
In this section, we consider a homogeneous space $\X=G/P$, where $G$ is a Lie group and $P$ a closed subgroup of $G$, such that the action of $P$ on $\X$ is faithfull.    On $\X$, we consider the canonical Cartan geometry defined as the triple $(\X,G,\omega^G)$, where $\omega^G$ is the Maurer-Cartan form on $G$.  A {\it Kleinian manifold} is a quotient $M=\Gamma \backslash \Omega$, where $\Gamma$ is a discrete subgroup of $G$ acting freely and properly on an open subset $\Omega \subset \X$.  We are thus in the situation of quotient manifolds described in section \ref{quotient-manifold}, where $(\X,G,\omega^G)$ plays the role of the Cartan geometry $(L,\hL,\ol)$.   If we refer to the notations introduced at the begining of section \ref{weak}, we will write $G$ (resp. $G_{\Gamma}=\Gamma \backslash G$) instead of $\hL$ (resp.  $\hL_{\Gamma}$). Also, we will adopt the notations $\omega^G$ and $\omega^{G_{\Gamma}}$ instead of $\ol$ and $\omega^{L_{\Gamma}}$ (the same for $d_L$, $\rho^L$, $\rho^{L_{\Gamma}}$ etc...). As explained in section \ref{weak}, there is a natural Cartan geometry $(M,\hm,\om)$ induced by $(\X,G,\omega^G)$, and when we will speak of {\it the} Cartan geometry  of a Kleinian manifold, we will always refer to this canonical one.

The fundamental fact here  is that the Riemannian manifold $(G,\rho^G)$ is $G$-homogeneous, hence {\it complete}.  Since $\pi_{\Gamma} : (G,\rho^{G}) \to (G_{\Gamma},\rho^{G_{\Gamma}}) $ is a Riemannian covering, the manifold $(G_{\Gamma},\rho^{G_{\Gamma}})$ is complete as well.  Thus, if $M=\Gamma \backslash \Omega$ is a Kleinian manifold and $(M,\hm,\om)$ its Cartan geometry,  $\hm$ is an open subset of the complete Riemannian manifold $(G_{\Gamma},\rho^{G_{\Gamma}})$, what will have interesting implications for geometrical embeddings of $(M,\hm,\om)$.

\subsection{The boundary map of a Kleinian manifold}
\label{application-bord-kleinienne}




As a first consequence of the completeness of $(G_{\Gamma},\rho^{G_{\Gamma}})$, we have a sharpening of lemma \ref{bord-normal}.

\begin{lemme}
\label{bord-kleinienne}
Let $\Omega  \subset \X$ be a normal domain of $\X$, and $M=\Gamma \backslash \Omega$  a  Kleinian manifold. Then the natural map $\lambda : \partial \hm \to \pchm$ is an homeomorphism. If $\sigma : (M,\hm,\om) \to (N,\hn,\on)$ is a strict geometrical embedding, then $\lambda$ maps $\hl$ onto $\hlc$.
\end{lemme}

\begin{preuve}
In Lemma \ref{bord-normal} we showed that $\lambda$ mapped $\partial \hm$ homeomorphically onto an open subset of $\partial \hm_c$. We are going  to prove that it is surjective. The key point here is that $(\gga,\rho^{\gga})$ is a complete Riemannian manifold. 
Let $\hx_{\infty}$ be a point of $\pchm$, and $(\hx_k)$ be a sequence of points of $\hm$ converging to $\hx_{\infty}$. The sequence $(\hx_k)$ is a Cauchy sequence for $d_M$, hence for $d_{G_{\Gamma}}$. By completeness of $(\gga,d_{\gga})$, $(\hx_k)$ converges to $\hx \in \phm$. By definition of $\lambda$, $\lambda(\hx)=\hx_{\infty}$, what proves that $\lambda$ is surjective. 

It is then straigthforward to check that if  $\sigma: (M,\hm\om) \to (N,\hn,\on)$ is a strict geometrical embedding, $\lambda$ maps $\hl$ onto $\hlc$.

\end{preuve}


\subsection{Proof of Theorem \ref{Cartan-codimension}}

We assume that $M= \Gamma \backslash \Omega$ is a Kleinian manifold of dimension $n$, modelled on $\X=G/P$, and that the Hausdorff dimension of $\partial \Omega$ is strictly less than $n-1$.  We consider $\sigma: (M,\hm,\omega^{M}) \to (N, \hn , \omega^N)$  a geometrical embedding.    We keep the notations of section \ref{application-bord-kleinienne}.  

The first observation is that because the Hausdorff dimension of $\hlc$ is smaller than $m-1$, where $m=dim(\hm)=dim(\lieg)$, and $\hl_c$ is open in $\partial \hm$, $\hm^{\prime}=\hlc \cup \hm$ is a dense open subset of $\gga$, on which $P$ acts freely and properly (proposition \ref{boundary}).  As we  saw in lemma \ref{identification}, $\Lambda_c=\Gamma \backslash \tilde \Lambda$, where $\tilde \Lambda$ is an open subset of $\partial \Omega$ on which $\Gamma$ acts freely properly discontinuously. Actually, $\Omega^{\prime}=\tilde \Lambda \cup \Omega$   
is also an open subset of $\X$ on which $\Gamma$ acts freely properly discontinuously (just because $P$ acts freely and properly on ${\tilde \Lambda}_c \cup \hm$). The space $M^{\prime}=\Lambda_c \cup M$ is thus a Kleinian manifold $M^{\prime}=\Gamma \backslash \Omega^{\prime}$, and if we denote by $\omega^{M^{\prime}}$ the restriction of $\omega^{\gga}$ to $\hm^{\prime}$, then $(M^{\prime},\hm^{\prime},\omega^{M^{\prime}})$ is the Cartan geometry of the Kleinian manifold $M^{\prime}$.  By the study made in section \ref{application-bord}, $\sigma$ extends to a continuous   map   $\overline \sigma : \hm^{\prime} \to \hn$.  For every $\hx \in \hl_c$, the open set $V_{\hx}$ of Proposition \ref{extension} can be choosen arbitrary small, so that in fact $V_{\hx} \subset \hm$. We thus see that on $V_{\hx}$, the map $\overline \sigma$ coincides with the map $\sigma_{\hx}^{\prime}$, which is smooth. We infer that $\overline \sigma$ is smooth, and ${\overline \sigma}^*\omega^N=\omega^{\hm^{\prime}}$, since this equality holds on $\hm$, which is dense in $\hm^{\prime}$. Finally, we have proved that $\overline \sigma : (M^{\prime},\hm^{\prime},\omega^{M^{\prime}}) \to (N,\hn,\omega^N)$ is a geometrical embeding. 

We would like to prove now that $\overline \sigma$ is a geometrical isomorphism.  By Proposition \ref{boundary}, it is sufficient to show that $\hl_c^{\prime}$ is empty (here the points of $\hl_c^{\prime}$ are those $\hx^{\prime} \in \partial \hm^{\prime}$ for which there exists a sequence $(x_k^{\prime})$ in $\hm^{\prime}$, such that $x_k^{\prime} \to \hx^{\prime}$ and $\overline \sigma(x_k^{\prime})$ converges in $\hn$).

Assume for a contradiction that there exists $\hx^{\prime} \in \hl_c^{\prime}$.  Let $(x_k^{\prime})$ be a sequence of $\hm^{\prime}$ converging to $\hx^{\prime}$, such that $\overline \sigma(\hx_k^{\prime})$ converges to $\hy \in \hn$.  Let $(x_k)$ a sequence of $\hm$
converging to $\hx^{\prime}$ (such a sequence exists since $\hm$ is dense in $\hm^{\prime}$). Since $\partial \Omega^{\prime} \subset \partial \Omega$, its Hausdorff dimension is $< n-1$.  It follows that $\Omega^{\prime}$ is  a normal domain of $\X$. Hence $\hx^{\prime}$ admits a neighbourhood $U \subset \gga$, such that $U \cap \hm^{\prime}$ is connected and  $d_{\gga} \leq d_{\gga}^{U \cap \hm^{\prime}} \leq k_{\hx^{\prime}} d_{\gga}$ on $U \cap \hm^{\prime}$, for a real $k_{\hx^{\prime}} \geq 1$.  As a consequence, $d_{M^{\prime}}^{U \cap \hm^{\prime}}(\hx_k^{\prime},\hx_k)= d_{\gga}^{U \cap \hm^{\prime}}(\hx_k^{\prime},\hx_k)$ converges to  $0$, what implies that $d_{\overline \sigma}^N(\overline \sigma(\hx_k^{\prime}),\overline \sigma(\hx_k))$ tends to $0$. Hence $\overline \sigma(\hx_k)=\sigma(\hx_k)$  tends to $\hy$, and $\hx^{\prime} \in \hl_c \subset \hm^{\prime}$, a contradiction.  

We get that $\overline \sigma : (M^{\prime},\hm^{\prime},\omega^{M^{\prime}}) \to (N,\hn,\omega^N)$ is a geometrical isomorphism, what proves the theorem.

\subsubsection{Remarks about theorem \ref{Cartan-codimension}}

Let us mention two interesting consequences of theorem \ref{Cartan-codimension}.  

Assume first that $\X$ is noncompact (this is the case if we are dealing with pseudo-Riemannian metrics of type $(p,q)$, in which case $\X$ is just the type-$(p,q)$ Minkowski space $\re^{p,q}$, or if we consider affine structures, for which $\X$ is the affine space).  Then given $\Omega \subset X$ an open subset such that $\partial \Omega$ has a Hausdorff dimension $< \dim X -1$. A straigthforward consequence of theorem \ref{Cartan-codimension} is that  {\it $\Omega$ does not admit any geometrical compactification}.  In other words, there is no geometrical embedding $\sigma : (\Omega, \hat \Omega, \omega^{\Omega}) \ to (N, \hn, \on)$, where $(N,\hn,\on)$ is a cartan geometry modelled on $\X$, and $N$ is compact.

On the other hand, let us assume that $\X$ is compact (for example if we are considering conformal pseudo-Riemannian structures of type $(p,q)$, or projective structures).  Then given $\Omega \subset X$ an open subset such that $\partial \Omega$ has a Hausdorff dimension $< \dim X -1$, theorem \ref{Cartan-codimension} implies {\it unicity of a geometrical compactification for $\Omega$}.  Namely, if $\sigma : (\Omega, \hat \Omega, \omega^{\Omega}) \to (N,\hn,\on)$ is a geometrical embedding, with $(N,\hn,\on)$ a Cartan geometry modelled on $\X$, $N$ compact, we get that $(N,\hn,\on)$ is geometrically isomorphic to $(\X,G,\omega^G)$.



\subsection{The case of domains bounded by hypersurfaces: Theorem \ref{klein-hypersurface}}

Here $\Omega \subsetneq \X$ is an open subset such that $\overline \Omega$ is a topological submanifold with boundary of $\X$, and $\partial \Omega$ is of class ${\cal C}^{k,\alpha}$, $k+\alpha \geq 1$. Let $M=\Gamma \backslash \Omega$ be a Kleinian manifold, and $\sigma: (M,\hm,\om)  \to (N,B,\omega^N)$ be a strict geometric embedding, where  $(N,B,\omega^N)$ is a Cartan geometry modelled on $\X$. We assume that the boundary $\partial_sM$ is locally the graph of a ${\cal C}^{0}$ map.  Lemma \ref{bord-kleinienne} ensures  that $\hl$ is a nonempty open subset of $\partial \hm$, hence a submanifold of $G_{\Gamma}$ of class ${\cal C}^{k,\alpha}$.  Since $\partial_sM$ is assumed to be locally a graph, every point of $\partial_sM$ is accessible, so that $\partial s(\Lambda)=\partial_sM$ by lemma \ref{accessible}. Let $y \in \partial_sM$, and $x \in \Lambda$ such that $\partial s(x)=y$.  Let $\hx \in \hat \Lambda$ and $\hy \in \hn$ above $x$ and $y$ respectively.  We apply proposition \ref{extension}.  There are $V_{\hx}$ and $V_{\hy}$ neighbourhoods of $\hx$ and $\hy$ in $G_{\Gamma}$ and $\hn$ respectively, and $\sigma_{\hx}^{\prime} : V_{\hx} \to V_{\hy}$ a smooth diffeomorphism inducing $\overline \sigma$ on $V_{\hx} \cap (\hl \cup \hm)$.  In particular, if ${\hat \Sigma}_{\hx}=V_{\hx} \cap \hl$, then $\partial \sigma({\hat \Sigma}_{\hx})$ is a ${\cal C}{k,\alpha}$-submanifold of $\hn$ containing $\hy$. By lemma \ref{identification}, there is  $\Sigma_x \subset \partial \Omega$ an open neighbourhood of $x$ (on which $\Gamma$ acts freely and properly), such that $\partial s(\Sigma_x)$ is the projection of $\partial \sigma({\hat \Sigma}_{\hx})$ on $N$, hence a ${\cal C}^{k, \alpha}$-submanifold of $N$, containing $y$.  Of course, $\partial s(\Sigma_x) \subset \partial_sM$, and since this latter is locally the graph of a map, we must conclude that $\partial s(\Sigma_x)$ is a neighbourhood of $y$ in $\partial_sM$, and that $\partial_sM$ is of class ${\cal C}^{k,\alpha}$.

\section{Maximal geometries: some examples}
\label{section-maximal}
\subsection{Proof of theorem \ref{klein-maximal}}

 Here $\Omega$ and $M=\Gamma \backslash \Omega$ satisfy the hypotheses of theorem \ref{klein-maximal}.  Assume for a contradiction that $\sigma : (M,\hm,\om) \to (N,\hn,\on)$ is a strict geometrical embedding.  Then, by lemma \ref{bord-kleinienne}, $\hl$ is a nonempty open subset of $\hl_c$.  Thus, lemma \ref{identification} ensures that there exists  $\tilde \Lambda$  a nonempty open subset of $\partial \Omega$, such that the action of $\Gamma$ on $\tilde \Lambda \cup \Omega$ is free and proper.  But by assumption, $\Gamma$ can not act freely and properly on somme open subset of $\overline \Omega$ containing strictly $\Omega$. This yields a contradiction and proves theorem \ref{klein-maximal}.

Let us also quote the:

\begin{proposition}
\label{complete-maximal}
Any Kleinian manifold of the form $M=\Gamma \backslash \X$ is geometrically maximal.
\end{proposition}   
\begin{preuve}
The proof is obvious with our previous work: the Cartan geometry $(M,\hm,\om)$ is {\it complete}, so that for any geometrical embedding $\sigma : (M,\hm,\om) \to (N,\hn,\on)$, $\hl_c$ must be empty.  By proposition \ref{boundary}, $\sigma$ can not be strict.
\end{preuve}

\subsection{Illustration of theorem \ref{klein-maximal}}
\label{applications}
The reader is referred to \cite{chern}, \cite{tanaka}, \cite{kobayashi}, \cite{sharpe} for the interpretation of conformal structures, $CR$-structures and projective structures as Cartan geometries.

\subsubsection{Proof of theorem \ref{courbure-constante-maximale}}

$\bullet$  {\it complete flat Riemannian structures}. A complete flat structure  is a quotient $\Gamma \backslash \Omega$, where $\Omega = {\bf S}^n \setminus \{ p \}$, and $\Gamma$ is a discrete subgroup of $PO(1,n+1)$ fixing $p$.  As soon as $\Gamma \not = \{ id \}$, the action of $\Gamma$ on $\partial \Omega$ is not free.  When $n \geq 3$, we can apply theorem \ref{klein-maximal} and deduce immediately point $1.$ of theorem \ref{courbure-constante-maximale}.

$\bullet$ {\it complete hyperbolic manifolds}.  Such manifolds are quotients $\Gamma \backslash \Omega$, where $\Omega$ is the upper-hemisphere in ${\bf S}^n$, and $\Gamma \subset PO(1,n+1)$ is a lattice.  If $\Lambda_{\Gamma}$ denotes the limite set of $\Gamma$ on $\partial \Omega$, it is well known that the action of $\Gamma$ is proper on ${\bf S}^n \setminus \Lambda_{\Gamma}$.  Hence, theorem \ref{klein-maximal} implies that $\Gamma \backslash {\bf H}^n$ is conformally maximal if and only if $\Lambda_{\Gamma}=\partial \Omega$.  This yields the first part of the second point of theorem \ref{courbure-constante-maximale}.
Recall now the projective model of the hyperbolic space ${\bf H}^n$.  Let $\re^{1,n}$ be the Minkowski space of dimension $n+1$, namely $\re^{n+1}$ endowed with the quadratic form $q^{1,n}(x)=-x_1^2+.....+x_{n+1}^2$.  Then the projective structure on ${\bf H}^n$ is the projective structure induced by that of $\bf{RP}^n$ on the projectivisation $\Omega$ of $\{  q^{1,n} < 0\}$.    The projectivization of the isotropic cone of $q^{1,n}$ is the topological boundary $\partial \Omega$, and is diffeomorphic to a $(n-1)$-sphere.  For any $x \in \partial \Omega$, one can define a hypersurface $\Sigma(x)$ as follows:  let $u_x$ be a vector of $\re^{1,n}$ projecting on $x$, then the $q^{1,n}$-orthogonal of $u_x$ is a degenerate hyperplane; one defines $\Sigma(x)$ as the projectivization of this hyperplane.  It is straigthforward to check that $\Sigma_x$ is tangent to $\partial \Omega$ at $x$, and that it divides  ${\bf RP}^n$ into two connected components, exactly one of each containing $\Omega$.  We call this component $\Omega(x)$.   The group $\Gamma$ is  a discrete subgroup of $PO(1,n) \subset PGL(n+1, \re)$. The condition $\Lambda_{\Gamma}={\bf S}^{n-1}$ implies that $\Gamma$ acts minimally on $\partial \Omega$.  It is then clear from theorem \ref{klein-maximal} that $\Gamma$ can not act properly on some open subset $\Omega^{\prime} \subset {\overline \Omega}$ containing strictly $\Omega$.  The hyperbolic manifold $\Gamma \backslash {\bf H}^n$ is then projectively maximal.

On the other hand, let us show that if $\Lambda_{\Gamma} \not = {\bf S}^{n-1}$, then the hyperbolic manifold $\Gamma \backslash {\bf H}^n$ is not projectively maximal.  Recall the decomposition $PO(1,n)=KAK$, where $K$ (resp. $A$) is a maximal compact subgroup (resp. a Cartan subgroup) of $PO(1,n)$.  Any sequence $(\gamma_k)$ in $\Gamma$ tending to infinity writes as $\gamma_k=l_{1,k}a_kl_{2,k}$, where $l_{1,k}$ and $l_{2,k}$ are two sequences of $K$ and $a_k$ tends to infinity in $A$.  One says that $(\gamma_k)$ tends simply to infinity when $l_{1,k}$ and $l_{2,k}$ both have a limit in $K$, that we call $l_{1,\infty}$ and $l_{2,\infty}$ respectively.  The group $A$ fixes exactly two points $p^+$ and $p^-$ on $\partial \Omega$, and acts freely properly on ${\bf RP}^n \setminus \{  \Sigma(p^+) \cup \{ p^- \}  \}$ and ${\bf RP}^n \setminus \{  \Sigma(p^-) \cup \{ p^+ \}  \}$. If $(\gamma_k)$ tends simply to infinity, we associate to it $p^+(\gamma_k)=l_{1,\infty}.p^+$ and $p^-(\gamma_k)=l_{2,\infty}^{-1}.p^-$. The action of $\Gamma$ is free and proper on $\Omega^{\prime}=\bigcup \{\Omega^+(p^+(\gamma_k))$,  the union being taken  over all sequences of $\Gamma$ tending simply to infinity.     Now, it is a classical fact that $\Lambda_{\Gamma}$ is actually the union of all the $p^+(\gamma_k)$, when  $(\gamma_k)$ ranges among all sequences of $\Gamma$ tending simply to infinity.  If $\Lambda_{\Gamma} \not = {\bf S}^{n-1}$, then $\Omega^{\prime}$ contains $\Omega$ as a strict open subset, and thus $\Gamma \backslash \Omega$ is not projectively maximal.

\subsubsection{Other applications}
We enumerate other examples illustrating theorem \ref{klein-maximal}.  

$\bullet$  {\it  Conformal Kleinian manifolds}.  Let us go further in the discussion above about hyperbolic manifolds. If $\Gamma$ is a noncompact discrete subgroup of $PO(1,n+1)$, then the action of $\Gamma$ on ${\bf S}^n$ splits into two parts.  There is a minimal closed $\Gamma$-invariant subset $\Lambda_{\Gamma}$ called the limit set of $\Gamma$.  On the complementary $\Omega_{\Gamma}$ (the so called domain of discontinuity), the action of $\Gamma$ is proper.  The group is {\it Kleinian} when $\Omega_{\Gamma} \not = \emptyset$.  A direct consequence of theorem \ref{klein-maximal} and the minimality of the action on $\Lambda_{\Gamma}$ is that Kleinian manifolds of the form $\Gamma \backslash \Omega_{\Gamma}$ are conformally maximal (among Riemannian conformal structures).

$\bullet$  {\it CR geometry}.  The boundary ${\bf S}^{2n-1}$ of the $n$-dimensional complex hyperbolic space ($n \geq 2$) is endowed with  a flat $CR$-structure, and actually, any $CR$-structure on a $(2n-1)$-dimensional manifold determines in a canonical way a Cartan geometry modelled on ${\bf S}^{2n-1}$.  The group $PU(1,n)$ is the group of $CR$-automorphisms of  ${\bf S}^{2n-1}$, and $Heis(2n-1) \subset PU(1,n)$ (the $(2n-1)$-dimensional Heisenberg group) acts simply transitively on the complementary of a point in  ${\bf S}^{2n-1}$.  Hence $Heis(2n-1)$ inherits also of a flat $CR$-structure, which is left-invariant.  A straigthforward consequence of theorem \ref{klein-maximal} is that {\it any $CR$-structure of the form $\Gamma \backslash Heis(2n-1)$, with $\Gamma \not = \{ id \}$ a discrete subgroup of $Heis(2n-1)$, is $CR$-maximal.}  

In the same way, the unit tangent bundle to the real hyperbolic space ${\bf H}^{n}$ ($n \geq 2$) is endowed with a flat $CR$-structure.  Actually, $T^1{\bf H}^n$ is an open orbit $\Omega$ of $O(1,n) \subset PU(1,n)$ in ${\bf S}^{2n-1}$, the boundary of which is a sphere ${\bf S}^{n-1}$.  The action of $O(1,n)$ on $\partial \Omega$ is identified to that on the boundary of ${\bf H}^n$.  As a consequence, we get that {\it the flat $CR$-structure on the unit tangent bundle of a complete hyperbolic manifold of finite volume is $CR$-maximal}.

$\bullet$  {\it Projectively maximal structures}.
The affine flat space $\re^n$ embeds projectively as an open subset $\Omega \subset {\bf RP}^n$.  The boundary $\partial \Omega$ is a codimension one ${\bf RP}^{n-1}$ in ${\bf RP}^n$.  Any translation of $\re^n$ extends as a projective transformation of ${\bf RP}^n$, which fixes pointwise $\partial \Omega$.  Applying theorem \ref{klein-maximal}, we get: {\it for $n \geq 3$, any complete flat affine manifold $\Gamma \backslash \re^n$, where $\Gamma$ is a nontrivial discrete subgroup of translations, is projectively maximal.}  

Also,  let us recall that by Bieberbach's theorem, a complete flat complete Riemannian manifold $\Gamma \backslash \re^n$, $\Gamma \not = \{ Id \}$, has a nontrivial subgroup of finite index which contains only pure translations.  Such a subgroup fixes pointwise $\partial \Omega$.  We thus get that except the $n$-dimensional Euclidean space $(n \geq 3)$, any flat complete Riemannian manifold is {\it projectively maximal}.



$\bullet$  {\it Examples of conformally maximal Lorentzian manifolds}. We conclude with a last example for  readers interested by Lorentzian geometry.  The {\it anti-de Sitter} space is the Lorentzian analogue of the real hyperbolic space.  A model for this space is the following: in $\re^{n+2}$ equiped with the quadratic form $q^{2,n}(x)=-x_1^2-x_2^2+....+x_{n+2}^2$, let us consider the quadric $Q_{-1}= \{  u \in \re^{n+2} \ | \ q^{2,n}(u)=-1 \}$.  The anti-de Sitter space is the quadric $Q_{-1}$ endowed with the Lorentz metric induced by $q^{2,n}$.  It is a complete Lorentz manifold of constant sectional curvature $-1$.   The group of isometries of ${\bf AdS}_{2n+1}$, $n \geq 1$, is $O(2,2n)$.  In particular, it contains subgroups isomorphic to $U(1,n)$.  Using the dynamical study made in section $3$ of \cite{frances2} (see also section $5$ of this reference), one can  prove that {\it if $\Gamma$ is a lattice in $U(1,n)$, then any anti-de Sitter manifold $\Gamma \backslash {\bf AdS}_{2n+1}$ is conformally maximal.}

\subsection{Maximal conformally homogeneous Riemannian structures: proof of Theorem \ref{homogene}}

We are considering $(M,g)$ a conformally homogeneous Riemannian manifold of dimension $\geq 3$, and we assume that $(M,g)$ is not conformally maximal, namely there exists a strict conformal embedding $s :(M,g) \to (N,h)$.  The proof will work in two steps.  The main step will be to prove the:

\begin{proposition}
\label{conf-plat}
If $(M,g)$ is conformally homogeneous but not conformally maximal, then it is conformally flat.  
\end{proposition}

Once this proposition will be proved, we will be reduced to the study of conformally flat conformally homogeneous Riemannian structures, which were investigated by D.V. Alekseevskii and B. N. Kimel'fel'd in \cite{aleks}. Since compact Riemannian structures are obviously conformally maximal, we restate the result they obtained for noncompact manifolds:

\begin{theorem} \cite{aleks}
\label{classification}
Let $(M,g)$ be a connected conformally flat Riemannain manifold of dimension $n \geq 3$, with a transitive group of conformal transformations, then $(M,g)$ is conformally eauivalent to one of the following spaces:  $\re^n$,${\bf H}^n$,${\bf H}^{n-1} \times \re$,${\bf H}^{n-m} \times {\bf S}^m$, $\re^{n-m} \times {\bf T}^m$, ${\bf H}^{n-m} \times Q^m$, where $Q^m$ is a nontrivial finite quotient of the round sphere ${\bf S}^m$. In all examples, $n-1 \geq m \geq 1$.
\end{theorem}

To deduce theorem \ref{homogene} from theorem \ref{classification}, we just have to rule out $\re^{n-m} \times {\bf T}^m$  and ${\bf H}^{n-m} \times Q^m$ ($n-1 \geq m \geq 1$), showing that there are conformally maximal.   The space $\re^{n-m} \times {\bf T}^m$ is just the quotient of $\re^n$ by a nontrivial discrete group of translations.  It is conformally maximal by the first point of theorem \ref{courbure-constante-maximale}.  The spaces ${\bf H}^{n-m} \times Q^m$ are quotients of  ${\bf H}^{n-m} \times {\bf S}^m$ by a group of the form $\{ id \} \times \Gamma \subset O(1,n-m) \times O(m+1)$, with $\Gamma$ finite and nontrivial in $O(m+1)$.  Now ${\bf H}^{n-m} \times {\bf S}^m$ is conformally equivalent to an open subset $\Omega$ of ${\bf S}^n$, obtained by removing a round sphere of dimension $(n-m-1)$ in ${\bf S}^n$ (with the convention that round spheres of dimension $0$ are the union of two points).  Hence  ${\bf H}^{n-m} \times Q^m$ are Kleinian manifolds.  Since it is easily checked that subgroups of the form $\{  id \} \times \Gamma$ as above fix $\partial \Omega$ pointwise, it is a consequence of theorem \ref{klein-maximal} that  ${\bf H}^{n-m} \times Q^m$ is conformally maximal when $\Gamma$ is nontrivial.

We now focus on the proof of proposition \ref{conf-plat}.  Actually, the proof is very much the same as that made in \cite{frances1} to generalize the Ferrand-Obata theorem.  We call $U$ the open set $s(M)$, so that the group $\Conf (M,g)$ is identified with the group of conformal transformations of $(U,h)$.  Let us call $(N,\hn,\on)$ the canonical Cartan geometry associated to the conformal structure $(N,h)$, and $\hat U$ the inverse image of $U$ by the projection $\pi_N: \hn \to N$. Every conformal transformation $f\in Conf(U)$ lifts to  $\hat f $, which acts by  bundle automorphism on $\hn$, s.t $(\hat f)^* \on_{|\hat U}=\on_{|\hat U}$. We assumed that $s$ is a strict conformal embedding, hence we can pick $y_{\infty} \in \partial U$.  Let us fix $x \in U$. By homogeneity of $U$, we can find $(f_k)$ a sequence of $\Conf U$  such that $y_k=f_k(x)$ tends to $y_{\infty}$. We choose $\hx$ and $\hy_{\infty}$ in $\hn$ projecting on $x$ and $y_{\infty}$ respectively. Then, there exists a sequence $(p_k)$ of $P$ such that $\hy_k={\hat f}_k(\hx).p_k^{-1}$ tends to $\hy$.  The sequence $(p_k)$ is called {\it an holonomy sequence} associated to ${\hat f}_k$.   The first property (\cite{frances1}, Theorem 3, point $(ii)$) is that since $(f_k)$ tends to infinity in $\Conf U$, then $(p_k)$ also tends to infinity.  Let $C_{\hx}({\cal B},\lambda)$ be a conformal cone included in $U$.  
We now apply Lemma \ref{dynamique-cones}:

$\bullet $ either there exists ${\cal B}^{\prime} \subset {\cal B}$ such that $p_k.C({\cal B}^{\prime},\lambda) \to o$.   The relation $ {\cal D}_{\hy_k}^{y_k}(C_{\hat f_k(\hx)}({\cal B},\lambda))=p_k. C({\cal B},\lambda)$ and Lemma \ref{to-zero} then imply that $f_k(C_{\hx}({\cal B},\lambda)) \to y_{\infty}$.    It is now quite a standard fact that the Weyl tensor (resp. the Cotton tensor if we are in dimension $3$) has to vanish on $C_{\hx}({\cal B},\lambda)$.     Indeed, $\int_{C_{\hx}({\cal B},\lambda)}|| W ||_h^{\frac{n}{2}}dVol_h$ is conformally invariant, where $W$ denotes de Weyl tensor of $(N,h)$. By the property $f_k(C_{\hx}({\cal B},\lambda)) \to y_{\infty}$, $\int_{C_{\hx}({\cal B},\lambda)}|| W ||_h^{\frac{n}{2}}dVol_h \leq \int_{B_h(y_{\infty}, \epsilon)}|| W ||_h^{\frac{n}{2}}dVol_h$, and this for every $\epsilon>0$ arbitrary small. This implies $\int_{C_{\hx}({\cal B},\lambda)}|| W ||_h^{\frac{n}{2}}dVol_h=0$, and finally $W=0$ on $C_{\hx}({\cal B},\lambda)$.   In dimension $3$, one considers $||C||$ instead of $||W||_h^{\frac{n}{2}}$.    Now, by homogeneity, we conclude that $W=0$ (resp.  $C=0$ if we are in dimension $3$) on $U$, what proves that $U$  is conformally flat.

$\bullet$  if we are in the second case of Lemma \ref{dynamique-cones}, there exist ${\cal C}({\cal B}_{\infty},\lambda_{\infty})$ a cone in $\lien^+$, $(l_k)$ a sequence of $P$ tending to $l_{\infty} \in P$, and $\epsilon_k \to 0$ such that:
$$ \Ad l_kp_k.{\cal C}({\cal B},\epsilon_k) \to {\cal C}({\cal B}_{\infty},\lambda_{\infty})$$
This implies, as in the proof of  Proposition \ref{property-H} that the $f_k(C_{\hx}({\cal B},\epsilon_k))$ contain, for $k$ large enough, $\overline V$, where  $V$ is a relatively compact open set of $N$.  We thus have $\overline V \subset f_k(C_{\hx}({\cal B},\epsilon_k)) \subset U$ for $k$ large, and since the sequence $(f_k)$ tends to infinity in $\Conf U$, this implies that the action of $\Conf U$ on $U$ is nonproper.  We can then use the  Ferrand-Obata theorem (\cite{ferrand}, \cite{obata}, \cite{schoen}, \cite{frances1}) which ensures that $(U,h)$ is conformally equivalent to the euclidean space.  In particular, it is conformally flat.

\subsection{About conformal compactifications: proof of Theorem \ref{gros-isometries}}
We consider a noncompact Riemannian manifold $(M,g)$ of dimension $n \geq 3$, having a noncompact group of isometries.  We assume that we have a (necessarily strict) conformal embedding $s : (M,g) \to (N,h)$, where $(N,h)$ is a compact $n$-dimensional Riemannian manifold.  As in the previous proof, we identify $(M,g)$ conformally with $(U,h)$, where $U=s(M)$.  Let us fix $x \in U$.  Since the isometry group of $(M,g)$ acts properly on $(M,g)$, the hypothesis that $\Is (M,g)$ is noncompact yields a sequence $(f_k)$ of $\Conf (U,h)$ tending to infinity such that $y_k=f_k(x)$ tends to $y_{\infty} \in y_{\infty}$.  We are now in the same situation as in the previous proof. Keeping the notations of this proof, we get that either $U$ is the euclidean space (hence is conformally flat), or the Weyl tensor (resp.  the Cotton tensor) has to vanish on some conformal cone $C_{\hx}({\cal B},\lambda)$ included in $U$.  In particular, this tensor vanishes at $x$ and since $x$ was any point of $U$, we conclude that $U$ is conformally flat.


\begin{remarque}
By the same techniques as those used in \cite{frances1}, one could prove analogue results for all rank one parabolic  geometry (for example $CR$-structures).  Namely, any homogeneous rank one parabolic geometry which is not maximal has to be locally flat.  Also, any rank one parabolic geometry with a noncompact automorphism group, which admits a geometric compactification (i.e a geometric embedding into some compact manifold), has to be locally flat.

\end{remarque}

\section{Annexe }
\label{annexe}
\subsection{Proof of lemma \ref{exemples-normaux}}
We assume first that the Hausdorff dimension of $\partial W$ is $< n-1$. Let us pick $x \in \partial W$, and consider a chart $\phi : V \to U$, where $V$ is a neighbourhood of $x$, and $U$ an open subset of $\re^n$ containing $0$.  We will assume that $\phi$ maps $x$ to $0$, and identify $ V$ with $U$ and $V \cap W$ with $\phi(V \cap W)$. Let us consider $B_r \subset V$ the ball of center $0$ and radius $r$, and let us take $x_0$ and $y_0$ in $B_r \cap W$. We denote by $W_{x_0}$ (resp. $W_{y_0}$) the  points in $B_r \cap W$ which can be joigned to $x_0$ (resp. $y_0$) by a line segment included in $B_r \cap W$. By lemma \ref{dense} (more exactly, by a straigthforward adaptation of its proof), $W_{x_0} \cap W_{y_0}$ is dense in $B_r \cap W$, so that we can find a sequence $\hz_k \in W_{x_0} \cap W_{y_0}$ which tends to $y_0$.  The line segment $[x_0,z_k] \subset B_r$ tends to the line segment $[x_0,y_0] \subset B_r$, and the line segment $[z_k,y_0] \subset B_r$ tends to a point. So, we get a sequence of broken line segments included in $B_r \cap W$, whose length tends to $d(x_0,y_0)$.  We thus have $d=d^{B_r}=d^{B_r \cap W}$ on $B_r \cap W$. We see that the neighbourhoods $U_i=B_{\frac{1}{i+1}}$ satisfy the definition \ref{normal-domains}, for $k_i=1$.  By remark \ref{remarque-normal}, the open subset $W$ is thus normal.

Assume now that $W$ is bounded by a  locally Lipschitz hypersurface. As above, looking at a suitable chart around $x \in \partial W$, we can restrict our study to a neighbourhood of $0$, $U=]-\epsilon, + \epsilon[ ^n \subset \re^n$. We see $\re^n$ as the product $\re^{n-1} \times \re$, and the boundary $W \cap U$ is given by a $k$-Lipschitz map $\phi : U^{n-1} \to \re$, where $U^{n-1}= ]-\epsilon, + \epsilon[ ^{n-1}$. Restricting to a suitably small neighbourhood of $x$, we will assume that the $(z,t) \in U \cap W$ are caracterized by $t > \phi(z)$, so that $U \cap W$ is connected. We endow $U$ with the Euclidean metric.  Since $\phi$ is $k$-Lipschitz, one has:  
$$d^{U^{\prime} \cap \hm} \leq \sqrt{1+k^2}d$$

The basis of neighbourhoods $]-\frac{\epsilon}{i+1}, + \frac{\epsilon}{i+1}[$, endowed with the Euclidean metric satisfy definition \ref{normal-domains}, if we put $k_i=k$.  By remark \ref{remarque-normal}, $W$ is a normal domain.

\subsection{Proof of lemma \ref{normal-fibre}}
Here we assume that $W \subset Y$ is a normal domain, $\pi: B \to Y$ is a fibration and $\hat W=\pi^{-1}(W)$.  We want to prove that $\hat W \subset B$ is also a normal domain.  Let $x \in W$, $\hx \in \hat W$ above $x$, and $(U_i)$ a sequence of neighbourhood as in definition \ref{normal-domains}, proving that $W$ is a normal domain of $Y$.  Since the property is local around $x$, there is no harm assuming that the $U_i$'s are open subsets of $\re^n$, and that $U_0$ is endowed with the Euclidean metric $g$ (see remark \ref{remarque-normal}). Also, in a suitable chart around $\hx$, $W_i=U_i \times ]-\epsilon_i,\epsilon_i[^m$ is a sequence of neighbourhoods of $\hx$ (with $m$ is the dimension of the fibers of $\pi$, and $\epsilon_i$  a decreasing sequence converging to $0$) satisfying $\bigcap_{i \geq 0} W_i = \{  \hx \}$,  and for all $i \geq 0$: $W_i \cap W$ is connected.  We endow $W_0$ with the product metric of the Euclidean metric   $g$ on $U_0$ and the Euclidean metric on $]-\epsilon_i,\epsilon_i[^m$. The resulting product metric is denoted by $\hat g$.  To avoid cumbersome notations, $L(\alpha)$ will indistinctively denote the length of $\alpha$ for $g$ if $\alpha$ is a curve of $U_i$, and for $\hat g$ if $\alpha$ is a curve of $W_i$.

Our aim is to prove that for every $i \geq 1$, there exists $K_i \geq 1$ such that $d^{W_i \cap W} \leq K_i d^{W_i}$.  this will be done thanks to the following lemmas.

\begin{lemme}
\label{l1}
Let $z_1=(x_1,y_1)$ and $z_2=(x_2,y_2)$ be two points of $W_i$  (with $x_j \in U_i$ and $y_j \in ]-\epsilon_i,\epsilon_i[^m$ ).
Let $\alpha = \alpha(u)$ be a ${\cal C}^1$-curve from $[0,L]$ to $U_i$ satisfying $g(\alpha^{\prime},\alpha^{\prime})=1$ and $\alpha(0)=x_1$, $\alpha(L)=x_2$. Then for every $\beta: [0,L] \to ]-\epsilon_i,\epsilon_i[^m$ of class ${\cal C}^1$ and satisfying $\beta(0)=y_1$ and $\beta(L)=y_2$, we have:
$$ L(\gamma) \geq \sqrt{L^2 + |y_2-y_1|^2}$$
where $\gamma: [0,L] \to W_i$ is defined by $\gamma(u)=(\alpha(u),\beta(u))$. The minimum $\sqrt{L^2 + |y_2-y_1|^2}$ is actually realized by $\beta(u)=(1-\frac{u}{L})y_1+\frac{u}{L}y_2$.

\end{lemme}  

\begin{preuve}
The length of $u \mapsto \gamma(u)=(\alpha(u),\beta(u))$ is just: $$L(\gamma)=\int_0^L\sqrt{1+|\beta^{\prime}(u)|^2}=L(\lambda_{\beta})$$ where $\lambda_{\beta} : [0,L] \to \re \times ]-\epsilon_i,\epsilon_i[^m$ joins the points $(0,y_1)$ and $(L,y_2)$ and is given by $\lambda_{\beta}(u)=(u,\beta(u))$.  It is clear that $L(\lambda_{\beta})$ achieves a minimum when  $\lambda_{\beta}$ is a straigth line, namely when $\beta(u)=(1-\frac{u}{L})y_1+\frac{u}{L}y_2$.  This minimal value is the Euclidean distance between $(0,y_1)$ and $(L,y_2)$, namely: $\sqrt{L^2+|y_1-y_2|^2}$.
\end{preuve}

We then prove:

\begin{lemme}
\label{l2}
Let $z_1=(x_1,y_1)$ and $z_2=(x_2,y_2)$ be two points of $W_i \cap \hat W$, $x_1 \not = x_2$.  Let $\gamma_k=(\gamma_{1,k},\gamma_{2,k})$  be a sequence of curves joining $z_1$ and $z_2$ in $W_i$ (resp. in $W_i \cap \hat W$), satisfying $\hat g (\gamma_{1,k}^{\prime},\gamma_{1,k}^{\prime})=1$, and such that $L_k=L(\gamma_k)$ tends to $d^{W_i}(z_1,z_2)$ (resp.  to $d^{W_i \cap \hat W}(z_1,z_2)$).  Then $L(\gamma_{1,k})$ tends to $d^{U_i}(x_1,x_2)$  (resp.  to $d^{U_i \cap W}$).
\end{lemme}

\begin{preuve}
Assume for a  contradiction that the conclusion does not hold.  Then we can find a curve $\alpha$ joining $x_1$ and $x_2$ in $U_i$  (resp. in $U_i \cap W$), and $\epsilon>0$ such that $L(\alpha)<L(\gamma_{1,k})-\epsilon$ for all $k \in \NN$.  By lemma \ref{l1}, there is a curve $\gamma$ joining $z_1$ and $z_2$ in $W_i$ (resp. in $W_i \cap \hat W$), projecting on $\alpha$, and having length $\sqrt{L(\alpha)^2+|y_1-y_2|^2}$.  On the other hand, the length of $\gamma_k$ is at least $\sqrt{L(\gamma_{1,k})^2+|y_1-y_2|^2}$.  We thus deduce the existence of $\delta>0$ such that $L(\gamma)<L(\gamma_k)-\delta$: contradiction. 

\end{preuve}

We now finish the proof of lemma \ref{normal-fibre}.  Let $z_1=(x_1,y_1)$ and $z_2=(x_2,y_2)$ two points of $W_i \cap W$.  If $x_1=x_2$, it is clear that $d^{W_1}(z_1,z_2)=d^{W_i \cap \hat W}(z_1,z_2)$.  So, we assume now that $x_1 \not = x_2$.  We choose $\gamma_k=(\gamma_{1,k},\gamma_{2,k})$  (resp. $\lambda_k=(\lambda_{1,k},\lambda_{2,k})$) a sequence of curves joining $z_1$ and $z_2$ in $W_i \cap \hat W$ (resp. in $W_i$) such that $g(\gamma_{1,k}^{\prime}, \gamma_{1,k}^{\prime})=1$  (resp.  $g(\lambda_{1,k}^{\prime}, \lambda_{1,k}^{\prime})=1$), and $L(\gamma_k) \to d^{W_i \cap \hat W}(z_1,z_2)$ (resp.  $L(\lambda_k) \to d^{W_i }(z_1,z_2)$).  By lemma \ref{l2}, we get that $L_{1,k}=L(\gamma_{1,k})$ tends to $d^{U_i \cap W}(x_1,x_2)$  (resp.  $l_{1,k}=L(\lambda_{1,k})$ tends to $d^{U_i }(x_1,x_2)$).  Moreover, by lemma \ref{l1}, we may also assume that $L(\gamma_k)=\sqrt{L_{1,k}^2+|y_1-y_2|^2}$  and $L(\lambda_k)=\sqrt{l_{1,k}^2+|y_1-y_2|^2}$. Since the $U_i$'s are as in  definition  \ref{normal-domains}, there is $k_i \geq 1$ such that $d^{U_i \cap W} \leq k_i d^{U_i}$.  Let us call $K_i=2k_i$.  Then for $k$ sufficiently large, $L_{1,k} \leq K_i l_{1,k}$, what yields:
$$ L(\gamma_k) \leq \sqrt{K_i^2l_{1,k}^2+|y_1-y_2|^2} \leq \sqrt{K_i^2(l_{1,k}^2+|y_1-y_2|^2)} \leq K_iL(\lambda_k)$$
Passing to the limit as $k \to +\infty$, we get:
$$ d^{W_i \cap \hat W}(z_1,z_2) \leq K_i d^{W_i}(z_1,z_2)$$
as desired.

\subsubsection{Acknowledgment}
This work was supported by the ANR {\it Geodycos}.

Charles FRANCES\\
Laboratoire de Math\'ematiques, Bat. 425.\\
Universit\'e Paris-Sud.\\
91405 ORSAY.\\
Charles.Frances@math.u-psud.fr

\end{document}